\documentclass[12pt]{article}%
\usepackage[dvips]{epsfig}
\usepackage{amsmath}
\usepackage{amsfonts}
\usepackage{amssymb}
\usepackage{graphicx}%
\numberwithin{equation}{section}

\renewcommand{\phi}{\varphi}

\newcommand{\Y}{\mathcal{Y}}

\renewcommand{\chi}{\mathcal{X}}

\begin{document}

\title{{\Large \textbf{Local estimates for elliptic equations arising in conformal
geometry }}}
\author{{\large Yan He\ \ \ \ \ \ Weimin Sheng }}
\date{\vspace{-8mm}}
\maketitle

\begin{abstract}
In this paper we consider Yamabe type problem for higher order curvatures on
manifolds with totally geodesic boundaries. We prove local gradient and second
derivative estimates for solutions to the fully nonlinear elliptic equations
associated with the problems.

\end{abstract}

\thispagestyle{empty}


\baselineskip=17.5pt \parskip=3pt

\section{Introduction}

\hskip600pt \footnote{The authors were supported by NSFC10771189 and 10831008.
The second author was also supported by the Australian Research Council.}

Let $(M^{n},g)$ be a smooth, compact Riemannian manifold of dimension $n\geq
3$. The Schouten tensor of $g$ is defined by
\[
A_{g}={\frac{1}{{n-2}}}\left(  \mathrm{Ric}_{g}-{\frac{{R_{g}}}{{2(n-1)}}%
}g\right)  ,
\]
where $\mathrm{Ric}_{g}$ and $R_{g}$ are the Ricci and scalar curvatures of
$g$, respectively. The $k$-curvature (or $\sigma_{k}$ curvature) is defined to
be the $k$-th elementary symmetric function $\sigma_{k}$ of the eigenvalues
$\lambda(g^{-1}A_{g})$ of $g^{-1}A_{g}$. If $\tilde{g}=e^{-2u}g$ is a metric
conformal to $g$, the Schouten tensor transforms according to the formula
\[
A_{\tilde{g}}=\nabla^{2}u+du\otimes du-{\frac{1}{2}}|\nabla u|^{2}g+A_{g},
\]
where $\nabla u$ and $\nabla^{2}u$ denote the gradient and Hessian of $u$ with
respect to $g$. Consequently, the problem of conformally deforming a given
metric to one with prescribed $\sigma_{k}$-curvature reduces to solving the
partial differential equation
\begin{equation}
\sigma_{k}\Big(\lambda\Big(g^{-1}\Big[\nabla^{2}u+du\otimes du-{\frac{1}{2}%
}|\nabla u|^{2}g+A_{g}\Big]\Big)\Big)=\psi(x)e^{-2ku}. \tag{1.1}\label{1.1}%
\end{equation}
For compact manifolds without boundary, the existence of the solutions to the
equation $\left(  \ref{1.1}\right)  $ has been studied by many authors (see
\cite{CGY1, CGY2, GW2, GW3, LL1, LL2, GV1, GV2, TW1, TW2, STW, GeW, V2} etc.)
since these equations were first introduced by J. A. Viaclovsky \cite{V1}.
$C^{1}$ and $C^{2}$ estimates have also been studied extensively, see
\cite{Cn1, GW1, GW2, LL1, STW, W2} for local interior estimates and \cite{V2}
for global estimates.

Another interesting problem is to study the fully nonlinear equation
(\ref{1.1}) on a compact Riemannian manifold $(M^{n},g)$ with boundary
$\partial M $. In \cite{G}, Bo Guan studied the existence problem under the
Dirichlet boundary condition. There are many poineering works on the Dirichlet
problems for fully nonlinear elliptic equations, see \cite{CNS, Tr2} etc.. The
Neumann problem for (\ref{1.1}) has been studied by S. Chen \cite{Cn2, Cn3},
Jin-Li-Li \cite{JLL}, Jin \cite{J} and Li-Li \cite{LL3}, etc.. Under various conditions,
they derive local estimates for solutions and establish some existence
results. Before introducing the problem, we need the following definitions.

Define
\[
\Gamma_{k}=\{\Lambda=(\lambda_{1},\lambda_{2},\dots,\lambda_{n})\in
\mathbb{R}^{n}|\quad\sigma_{j}\left(  \Lambda\right)  >0,1\leq j\leq k\},
\]
and $1\leq k\leq n$, where $\sigma_{k}$ is the $k$-th elementary symmetric
function defined by
\[
\sigma_{k}(\Lambda)=\sum_{1\leq i_{1}<\cdots<i_{k}\leq n}\lambda_{i_{1}}%
\cdots\lambda_{i_{k}}%
\]
for all $\Lambda=(\lambda_{1},\lambda_{2},\dots,\lambda_{n})\in\mathbb{R}^{n}$
. We also denote $\sigma_{0}=1$. Therefore we have the relation $\Gamma
_{n}\subset\Gamma_{n-1}\subset\cdot\cdot\cdot\subset\Gamma_{1} $. For a
$2$-symmetric form $S$ defined on $(M^{n},g)$, $S\in\Gamma_{k}$ means that the
eigenvalues of $S$, $\lambda\left(  g^{-1}S\right)  $ lie in $\Gamma_{k}$. We
also denote $\Gamma_{k}^{-}=-\Gamma_{k}$.

Let $(M^{n},g),n\geq3$, be a smooth compact Riemannian manifold with nonempty
smooth boundary $\partial M$. We denote the mean curvature and the second
fundamental form of $\partial M$ by $h_{g}$ and $L_{\alpha\beta}$, where
$\{x^{\alpha}\}_{1\leq\alpha\leq n-1}$ is the local coordinates on the
boundary $\partial M$, and $\frac{\partial}{\partial x^{n}}$ the unit inner
normal with respect to the metric. In this paper similar as \cite{E2, Cn2} we
use Fermi coordinates in a boundary neighborhood. In these local coordinates,
we take the geodesic in the inner normal direction $\nu=\frac{\partial
}{\partial x^{n}}$ parameterized by arc length, and $\left(  x^{1}%
,...,x^{n-1}\right)  $ forms a local chart on the boundary. The metric can be
expressed as $g=g_{\alpha\beta}dx^{\alpha}dx^{\beta}+\left(  dx^{n}\right)
^{2}$. The Greek letters $\alpha,\beta,\gamma,...$ stand for the tangential
direction indices, $1\leq\alpha,\beta,\gamma,...\leq n-1$, while the Latin
letters $i,j,k,...$stand for the full indices, $1\leq i,j,k,...\leq n$. In
Fermi coordinates, the half ball is defined by $\overline{B}_{r}^{+}%
=\{x_{n}\geq0,\sum_{i}x_{i}^{2}\leq r^{2}\}$ and the segment on the boundary
by $\Sigma_{r}=\{x_{n}=0,\sum_{i}x_{i}^{2}\leq r^{2}\}$. Under the conformal
change of the metric $\widetilde{g}=e^{-2u}g$, the second fundamental form
satisfies
\[
\widetilde{L}_{\alpha\beta}e^{u}=\frac{\partial u}{\partial\nu}g_{\alpha\beta
}+L_{\alpha\beta}.
\]
The boundary is called umbilic if the second fundamental form $L_{\alpha\beta
}=\tau_{g}g_{\alpha\beta}$, where $\tau_{g}$ is the function defined on
$\partial M$. A totally geodesic boundary is umbilic with $\tau_{g}\equiv0$.
Note that the umbilicity is conformally invariant. When the boundary is
umbilic, the above formula becomes
\[
\tau_{\widetilde{g}}e^{-u}=\frac{\partial u}{\partial\nu}+\tau_{g}.
\]
The $k$-Yamabe problem with umbilic boundary becomes to considering the
following equation:
\begin{equation}
\left\{
\begin{array}
[c]{lr}%
\sigma_{k}^{1/k}\Big(\lambda\Big(g^{-1}\Big[\nabla^{2}u+du\otimes du-{\frac
{1}{2}}|\nabla u|^{2}g+A_{g}\Big]\Big)\Big)=e^{-2u}\text{ } & \text{in }M,\\
\frac{\partial u}{\partial\nu}=\tau_{\widetilde{g}}e^{-u}-\tau_{g}\text{
\ \ \ \ \ \ \ \ \ \ \ \ \ \ \ \ \ \ \ \ \ \ \ \ \ \ \ \ \ \ \ \ \ \ \ \ \ \ \ \ \ \ \ \ \ \ \ \ \ \ \ \ \ \ \ }
& \text{on }\partial M.
\end{array}
\right.  \tag{1.2}\label{1.2}%
\end{equation}
In \cite{Cn2,Cn3}, \cite{JLL} and \cite{J}, the authors established the a priori
estimates and obtained some existence results for (\ref{1.2}).

In this paper, we will generalize their results to more general equations,
which in particular include the equation (\ref{1.2}). In \cite{GV3}, Gursky
and Viaclovsky introduced a modified Schouten tensor
\[
A_{g}^{t}={\frac{1}{{n-2}}}\left(  \mathrm{Ric}_{g}-{\frac{t{R_{g}}}{{2(n-1)}%
}}g\right)  ,
\]
where $t\in\mathbb{R}$ is a parameter. When $t=1$, $A_{g}^{1}$ is just the
Schouten tensor; $t=n-1$, $A_{g}^{n-1}$ is the Einstein tensor; while $t=0$,
$A_{g}^{0}$ is the Ricci tensor. This tensor $A_{g}^{t}$ is in fact a constant
multiple of the tensor $sA_{g}+\frac{(1-s)}{2(n-1)}R_{g}g$ which is introduced
in \cite{LL1}, i.e. $A_{g}^{\frac{2s-1}{s}}=\frac{1}{s}(sA_{g}+\frac
{(1-s)}{2(n-1)}R_{g}g)$. Under the conformal change of the metric
$\widetilde{g}=e^{-2u}g$, $A_{\widetilde{g}}^{t}$ satisfies
\[
A_{\widetilde{g}}^{t}=A_{g}^{t}+\nabla^{2}u+\frac{1-t}{n-2}(\triangle
u)g+du\otimes du-\frac{2-t}{2}|\nabla u|^{2}g.
\]
In \cite{LS} and \cite{SZ}, we have studied
\[
\sigma_{k}\left(  \lambda\left(  g^{-1}\left[  A_{g}^{t}+\nabla^{2}%
u+\frac{1-t}{n-2}(\triangle u)g+du\otimes du-\frac{2-t}{2}|\nabla
u|^{2}g\right]  \right)  \right)  =f\left(  x\right)  e^{-2ku}%
\]
for $t\leq1$ or $t\geq n-1$. By use of the parabolic approach, we obtained
some existence results. Let $(M,g)$ be a compact, connected Riemannian
manifold of dimension $n\geq3$ with umbilic boundary $\partial M$, $W$ be a
$\left(  0,2\right)  $ symmetric tensor on $\left(  M^{n},g\right)  $.
Motivated by \cite{Cn1}, in this paper we study the following equation
\begin{equation}
\left\{
\begin{array}
[c]{ll}%
F(g^{-1}W)=f(x,u) & \text{ \ in \ }M\text{\ }\\
\frac{\partial u}{\partial\nu}=\widetilde{\tau}e^{-u}-\tau & \text{ \ on
\ }\partial M
\end{array}
\right.  \tag{1.3}\label{1.3}%
\end{equation}
where $F$ satisfies some fundamental structure conditions listed below, and
$\tau$ is the principal curvature of the boundary $\partial M$. We will
establish local a priori estimates for the solutions to the equation
(\ref{1.3}). After that, we will give some applications. More applications,
see \cite{HS1, HS2}.

We now describe the fundamental structure conditions for $F$.

Let $\Gamma$ be an open convex cone with vertex at the origin satisfying
$\Gamma_{n}\subset\Gamma\subset\Gamma_{1}$. Suppose that $F(\lambda)$ is a
homogeneous symmetric function of degree one in $\Gamma$ normalized with
$F(e)=F((1,\cdots,1))=1.$ Moreover, $F$ satisfies the following in $\Gamma$:

(A1) $F$ is positive.

(A2) $F$ is concave (i.e., $\frac{\partial^{2}F}{\partial\lambda_{i}%
\partial\lambda_{j}}$ is negative semi definite).

(A3) $F$ is monotone (i.e., $\frac{\partial F}{\partial\lambda_{i}}$ is positive).

(A4) $\frac{\partial F}{\partial\lambda_{i}}\geq\varepsilon\frac{F}{\sigma
_{1}},$ for some constant $\varepsilon>0,$ for all $i$.

The conditions (A1), (A2), (A3) and (A4) are similar as those in \cite{Cn2}.
Before stating the theorems, we introduce the following notations. Let
$f(x,z):M^{n}\times\mathbb{R}\rightarrow\mathbb{R}$ be a given positive
function. Let $u=u(x):M^{n}\rightarrow\mathbb{R}$ be a solution to
(\ref{1.3}). We define
\[
\overline{c_{sup}}(r)=sup_{\overline{B}_{r}^{+}}(f+|\nabla_{x}f(x,u)|+|f_{z}%
(x,u)|+|\nabla_{x}^{2}f(x,u)|+|\nabla_{x}f_{z}(x,u)|+|f_{zz}(x,u)|)
\]
or%
\[
c_{sup}(r)=sup_{B_{r}}(f+|\nabla_{x}f(x,u)|+|f_{z}(x,u)|+|\nabla_{x}%
^{2}f(x,u)|+|\nabla_{x}f_{z}(x,u)|+|f_{zz}(x,u)|),
\]
which varies with boundary or interior estimates.

Now we turn to the first equation: let
\begin{equation}
W=\nabla^{2}u+\frac{1-t}{n-2}(\triangle u)\ g+a(x)du\otimes du+b(x)|\nabla
u|^{2}g+S, \tag{1.4}\label{1.4}%
\end{equation}
where $t$ is a constant satisfying $t\leq1$, $S$ a $2$-symmetric form defined
on $M$, and $a\left(  x\right)  $, $b\left(  x\right)  $ are two smooth
functions on $M$. The derivatives are covariant with respect to the metric
$g$. We have \\[0.2cm]
\textbf{Theorem 1.} \textit{Let }$F$\textit{\ satisfy
the structure conditions (A1)-(A4) in a corresponding cone }$\Gamma$,
$a\left(  x\right)  =a,b\left(  x\right)  =b$ \textit{are two constants,
}$S=A$ \textit{the Schouten tensor. Suppose that the boundary }$\partial
M$\textit{\ is totally geodesic. Let }$u(x)$\textit{\ be a }$C^{4}%
$\textit{\ solution to the equation }%
\begin{equation}
\left\{
\begin{array}
[c]{ll}%
F(g^{-1}(\nabla^{2}u+\frac{1-t}{n-2}\triangle u\ g+a(x)du\otimes
du+b(x)|\nabla u|^{2}g+S))=f(x,u) & \text{ in\ }\overline{B}_{r}^{+},\\
\frac{\partial u}{\partial x^{n}}=0 & \text{ on }\Sigma_{r},
\end{array}
\right.  \tag{1.5}\label{1.5}%
\end{equation}
\textit{\ and} $W\in \Gamma$\textit{. Suppose that} $|\nabla f|\leq\Lambda f,\ |f_{z}|\leq\Lambda
f$\textit{\ for some constant} $\Lambda>0$\textit{.} \textit{If} $\frac
{1-t}{n-2}a-b\ge\delta_{1}>0$, $a+nb\le-\delta_{3}<0$ \textit{, and} $a\geq
0$\textit{, then}
\[
sup_{\overline{B}_{\frac{r}{2}}^{+}}\left(  |\nabla^{2}u|+|\nabla
u|^{2}\right)  \leq C,
\]
\textit{where }$C$\textit{\ depends on} $r,n,\varepsilon,\Lambda,\delta
_{1},\delta_{3},a$, $b$, $||A||_{C^{2}(\overline{B}_{r}^{+})}$, $||g||_{C^{3}%
(\overline{B}_{r}^{+})}$\textit{\ and }$\overline{c_{sup}}(r)$\textit{.}
\newline

When $t=1$, and $a=1$, $b=-\frac{1}{2}$, the boundary estimates have been
obtained by S. Chen \cite{Cn2, Cn3}, Jin-Li-Li \cite{JLL} and Jin \cite{J} for some special
cases. When $t=1$, the local interior estimates have been discussed by S. Chen
in \cite{Cn1} for general functions $a\left(  x\right)  $, $b\left(  x\right)
$ and a general $2$-symmetric tensor $S$. We just focus on the interior
estimates for the same equation, we may get \\[0.2cm]
\textbf{Theorem 2. }
\textit{Let }$F$\textit{\ satisfy the structure conditions (A1)-(A4) in a
corresponding cone }$\Gamma$\textit{. Let }$u(x)$\textit{\ be a }$C^{4}%
$\textit{\ solution to the equation }%
\begin{equation}
F(g^{-1}(\nabla^{2}u+\frac{1-t}{n-2}\triangle u\ g+a(x)du\otimes
du+b(x)|\nabla u|^{2}g+S))=f(x,u)
\tag{1.6}\label{1.6}
\end{equation}
\textit{\ in a local geodesic ball }$B_{r}\subset M$\textit{\ and} $W\in \Gamma$\textit{. Suppose that }$|\nabla
f|\leq\Lambda f,\ |f_{z}|\leq\Lambda f$\textit{\ for some constant }${\Lambda
}>0$. \newline Case (a). \textit{If }$\Gamma\subset\Gamma_{2}^{+}$\textit{,
}$\frac{1-t}{n-2}a(x)-b(x)\ge\delta_{1}>0$\textit{\ and }$\min\left\{
2ab+b^2,b^{2}\right\}  \geq\delta_{2}>0,$\textit{then}
\[
sup_{{B}_{\frac{r}{2}}}\left(  |\nabla^{2}u|+|\nabla u|^{2}\right)  \leq C,
\]
\textit{where }$C$\textit{\ depends only on }$r,n,\Lambda,\delta_{1}%
,\delta_{2},||a||_{C^{2}({B}_{r})}$, $||b||_{C^{2}({B}_{r})}$, $||S||_{C^{2}%
({B}_{r})}$, $||g||_{C^{3}({B}_{r})}$\textit{\ and }${c_{sup}}(r)$.\newline
Case (b). \textit{If }$\frac{1-t}{n-2}a(x)-b(x)\ge\delta_{1}>0$\textit{,
}$a(x)+nb(x)\le-\delta_{3}<0$\textit{ and }$a\left(  x\right)  \geq0$,\textit{
then we have}
\[
sup_{B_{\frac{r}{2}}}\left(  |\nabla^{2}u|+|\nabla u|^{2}\right)  \leq C,
\]
\textit{where C depends on }$r,n,\varepsilon,\Lambda,\delta_{1},\delta
_{3},||a||_{C^{2}(B_{r})}$, $||b||_{C^{2}(B_{r})}$, $||S||_{C^{2}(B_{r})}$,
$||g||_{C^{3}(B_{r})}$ \textit{and }$c_{sup}(r)$.\\[0.2cm]
\textit{Remark.}\ In case (a), the condition $\min \left\{2ab+b^2, b^2\right\}\ge \delta_2>0$ may be replaced by $\min\left\{b^2+2ab-2n||a||_{\infty}b-a^2, b^2\right\}\ge \delta_2>0$. The proof can be found in the proof of Theorem 2, case (a). The later condition is better than the former when $a>0$.

The a priori estimates in Theorem 1 and Theorem 2 rely on the signs of
$a\left(  x\right)  $ and $b\left(  x\right)  $. In fact, in \cite{STW} the
authors give a counterexample to show that there is no regularity if $a(x)=0$
and $b\left(  x\right)  >0$ when $t=1$. It is well known that the equation
(\ref{1.6}) has another elliptic branch, namely when the eigenvalues $\lambda$
lie in the negative cone $\Gamma_{k}^{-}$. Now we consider the second
equation. Let
\[
V=\frac{t-1}{n-2}(\triangle u)\ g-\nabla^{2}u-a(x)du\otimes du-b(x)|\nabla
u|^{2}g+S,
\]
where $t$ is a constant satisfying $t\geq n-1$. We have \\[0.2cm]
\textbf{Theorem
3. }\textit{\ Let }$F$\textit{\ satisfy the structure conditions (A1)-(A4),
}$a\left(  x\right)  =a,b\left(  x\right)  =b$ \textit{are two constants,
}$-S=A$ \textit{the Schouten tensor. Suppose that the boundary }$\partial
M$\textit{\ is totally geodesic. Let }$u(x)$\textit{\ be a }$C^{4}%
$\textit{\ solution to the equation }%
\begin{equation}
\left\{
\begin{array}
[c]{ll}%
F(g^{-1}(\frac{t-1}{n-2}(\triangle u)g-\nabla^{2}u-a(x)du\otimes
du-b(x)|\nabla u|^{2}g+S))=f(x,u) & \text{ in\ }\overline{B}_{r}^{+},\\
\frac{\partial u}{\partial x^{n}}=0 & \text{on }\Sigma_{r},
\end{array}
\right.  \tag{1.7}\label{1.7}%
\end{equation}
\textit{and $V\in \Gamma$, $t> n-1$. Suppose that }$|\nabla f|\leq\Lambda f,\ |f_{z}|\leq\Lambda
f$\textit{\ for some constant }$\Lambda>0$. \textit{If } $a+nb\geq\delta
_{3}>0,a\geq0$\textit{,\ then}
\[
sup_{\overline{B}_{\frac{r}{2}}^{+}}\left(  |\nabla^{2}u|+|\nabla
u|^{2}\right)  \leq C,
\]
\textit{where }$C$\textit{\ depends on }$r,n,\varepsilon,\Lambda, \delta
_{3},a$, $b$, $||A||_{C^{2}(\overline{B}_{r}^{+})}$, $||g||_{C^{3}%
(\overline{B}_{r}^{+})}$ \textit{and }$\overline{c_{sup}}(r)$.

Similar with Theorem 2, if we just focus on the interior estimates for the
same equation, we can get the following theorem for general functions
$a\left(  x\right)  $, $b\left(  x\right)  $ and general $2$-symmetric tensor
$S$.
\newline\textbf{Theorem 4. } \textit{Let }$F$\textit{\ satisfy the
structure conditions (A1)-(A4) in a corresponding cone }$\Gamma$\textit{. Let
}$u(x)$\textit{\ be a }$C^{4}$\textit{\ solution to the equation }%
\begin{equation}
F(g^{-1}(\frac{t-1}{n-2}\left(  \triangle u\right)  g-\nabla^{2}%
u-a(x)du\otimes du-b(x)|\nabla u|^{2}g+S))=f(x,u)
\tag{1.8}\label{1.8}%
\end{equation}
\textit{in a local geodesic ball }$B_{r}\subset M$\textit{\ and $V\in \Gamma$. Suppose that }$|\nabla
f|\leq\Lambda f,\ |f_{z}|\leq\Lambda f$\textit{\ for some constant }%
$\Lambda>0$.
\newline Case (a).\textit{\ If }$\Gamma\subset\Gamma_{2}^{+}%
$\textit{, }$\frac{t-1}{n-2}a(x)+b(x)\geq\delta_{1}>0$\textit{, }$\min\left\{
2ab+b^{2},b^{2}\right\}  \geq\delta_{2}>0,$ \textit{then}
\[
sup_{{B}_{\frac{r}{2}}}\left(  |\nabla^{2}u|+|\nabla u|^{2}\right)  \leq C,
\]
\textit{where }$C$\textit{\ depends only on }$r,n,\Lambda,\delta_{1}%
,\delta_{2},||a||_{C^{2}({B}_{r})}$, $||b||_{C^{2}({B}_{r})}$, $||S||_{C^{2}%
({B}_{r})}$, $||g||_{C^{3}({B}_{r})}$\textit{\ and }${c_{sup}}(r)$\textit{.}
\newline Case (b).\textit{\ If } $a(x)+nb(x)\geq\delta_{3}>0, a(x)\geq0,$
\textit{then}
\[
sup_{{B}_{\frac{r}{2}}}\left(  |\nabla^{2}u|+|\nabla u|^{2}\right)  \leq C,
\]
\textit{where }$C$\textit{\ depends on} $r,n,\Lambda, \delta_{3}%
,||a||_{C^{2}({B}_{r})}$, $||b||_{C^{2}({B}_{r})}$, $||S||_{C^{2}({B}_{r})}$,
$||g||_{C^{3}({B}_{r})}$\textit{\ and }${c_{sup}}(r)$.

Our idea of proof is from \cite{Cn1, Cn2, Cn3}, that is we estimate the
quantity $K:=\Delta u+a(x)|\nabla u|^{2}$ rather than estimate the gradient
and second derivatives separately. This idea was first used by Sophie Chen in
\cite{Cn1}. As in \cite{Cn2, Cn3} we show that the function $Ke^{px_{n}}$ does
not attain its maximum on the boundary, where $x_{n}$ is the distance to the
boundary. The main point in our argument is the observation that there exists
a suitable conformal transformation such that the metric has some nice
geometric properties on the boundary (Lemma 5). We would like to mention a
different method in getting the boundary estimates \cite{JLL} and \cite{J}. For the Neumann
problem of the Monge-Amp\`ere equation, the estimates were first obtained in
\cite{LTU}.

This paper is organized as follows. We begin with some background in Section
2. In Section 3, we discuss the applications which are based on the a priori
estimates in Theorem 1 to get the existence result of $k$-Yamabe problem. In
Section 4 and Section 5 we first prove Theorem 2 and Theorem 4 respectively.
We then prove the maximum of $K$ does not appear on the boundary, therefore
Theorem 1 and Theorem 3 can be concluded by the similar arguments of the case
(b) of Theorems 2 and 4 respectively.

\section{Preliminaries}

In this section, we give some basic facts about homogeneous symmetric
functions and show some outcomes by direct calculation under Fermi
coordinates. All of the facts can be found in the literatures cited below.

From Lemma 1 and Lemma 2 below, we can conclude that $F$ satisfies
(A1)-(A4).\\[0.2cm]\textbf{Lemma 1. }\textit{(\cite{U}) Let }$\Gamma
$\textit{\ be an open convex cone with vertex at the origin satisfying
}$\Gamma_{n}^{+}\subset\Gamma$\textit{, and let }$e=(1,\cdots,1)$\textit{\ be
the identity. Suppose that }$F$\textit{\ is a homogeneous symmetric function
of degree one normalized with }$F(e)=1$\textit{, and that }$F$\textit{\ is
concave in }$\Gamma.$\textit{\ Then }\newline\textit{(a) }$\sum_{i}\lambda
_{i}\frac{\partial F}{\partial\lambda_{i}}=F(\lambda),$\textit{\ for }%
$\lambda\in\Gamma;$\newline\textit{(b) }$\sum_{i}\frac{\partial F}%
{\partial\lambda_{i}}\geq F(e)=1$\textit{, for }$\lambda\in\Gamma.$
\\[0.2cm]\textbf{Lemma 2. }\textit{(\cite{Tr2, LT})\ Let }$G=\left(
\frac{\sigma_{k}}{\sigma_{l}}\right)  ^{\frac{1}{k-l}},0\leq l<k\leq
n.$\textit{\ Then \newline(a) }$G$\textit{\ is positive and concave in
}$\Gamma_{k};$\textit{\newline(b) }$G$\textit{\ is monotone in }$\Gamma_{k}%
,$\textit{\ i.e., the matrix }$G^{ij}=\frac{\partial G}{\partial W_{ij}}%
$\textit{\ is positive definite;\newline(c) Suppose }$\lambda\in\Gamma_{k}%
.$\textit{\ For }$0\leq l<k\leq n,$\textit{\ the following is the
Newton-Maclaurin inequality }%
\[
k(n-l+1)\sigma_{l-1}\sigma_{k}\leq l(n-k+1)\sigma_{l}\sigma_{k-1}.
\]

The following two lemmas will be used in proving Theorem 1 and 2. Let us
review some formulas on the boundary under Fermi coordinates (see \cite{E2} or
\cite{Cn3}). The metric is expressed as $g=g_{\alpha\beta}dx^{\alpha}%
dx^{\beta}+\left(  dx^{n}\right)  ^{2}$. The Christoffel symbols satisfy
\[
\Gamma_{\alpha\beta}^{n}=L_{\alpha\beta},\Gamma_{\alpha n}^{\beta}%
=-L_{\alpha\gamma}g^{\gamma\beta},\Gamma_{\alpha n}^{n}=0,\Gamma_{\alpha\beta
}^{\gamma}=\widetilde{\Gamma}_{\alpha\beta}^{\gamma},
\]
on the boundary, where we denote the tensors and covariant derivations with
respect to the induced metric on the boundary by a tilde (e.g. $\widetilde
{\Gamma}_{\alpha\beta}^{\gamma}$, $\tau_{\tilde{\alpha}\tilde{\beta}}$). When
the boundary is umbilic, we have
\[
\Gamma_{\alpha\beta}^{n}=\tau g_{\alpha\beta},\Gamma_{\alpha n}^{\beta}%
=-\tau\delta_{\alpha\beta},\Gamma_{\alpha n}^{n}=0.
\]
\\[0.2cm]\textbf{Lemma 3. }\textit{(see \cite{Cn3}) Suppose boundary
}$\partial M$\textit{\ is umbilic. Let }$u$\textit{\ satisfy }$u_{n}%
:=\frac{\partial u}{\partial x^{n}}=-\tau+\widetilde{\tau}e^{-u}%
,$\textit{\ where }$\widetilde{\tau}$\textit{\ is constant. Then on the
boundary we have }
\begin{equation}
u_{n\alpha}=-\tau_{\alpha}+\tau u_{\alpha}-\widetilde{\tau}u_{\alpha}e^{-u};
\tag{2.1}\label{2.1}%
\end{equation}
\textit{and}
\begin{align}
u_{\alpha\beta n}=  &  (2\tau-\widetilde{\tau}e^{-u})u_{\alpha\beta}-\tau
u_{nn}g_{\alpha\beta}+\widetilde{\tau}u_{\alpha}u_{\beta}e^{-u}-\tau
_{\tilde{\alpha}\tilde{\beta}}+\tau_{\alpha}u_{\beta}+\tau_{\beta}u_{\alpha
}\nonumber\\
&  -\tau_{\gamma}u_{\gamma}g_{\alpha\beta}+R_{n\beta\alpha n}(-\tau
+\widetilde{\tau}e^{-u})-\tau(-\tau+\widetilde{\tau}e^{-u})^{2}g_{\alpha\beta
}. \tag{2.2}\label{2.2}%
\end{align}
\\[0.2cm]\textbf{Lemma 4. }\textit{Suppose the boundary }$\partial
M$\textit{\ is totally geodesic and }$u_{n}=0$\textit{\ on the boundary. Then
we have}{\ }\textit{on the boundary}%
\begin{equation}
W_{\alpha\beta n}=\frac{1-t}{n-2}u_{nnn}g_{\alpha\beta}+a_{n}u_{\alpha
}u_{\beta}+b_{n}\left(  \Sigma_{\gamma}u_{\gamma}^{2}\right)  g_{\alpha\beta
}+S_{\alpha\beta n}, \tag{2.3}\label{2.3}%
\end{equation}%
\begin{equation}
V_{\alpha\beta n}=\frac{t-1}{n-2}u_{nnn}g_{\alpha\beta}-a_{n}u_{\alpha
}u_{\beta}-b_{n}\left(  \Sigma_{\gamma}u_{\gamma}^{2}\right)  g_{\alpha\beta
}+S_{\alpha\beta n}. \tag{2.4}\label{2.4}%
\end{equation}
\newline\textit{Proof.}\textbf{ }By the boundary condition we know that
$\widetilde{\tau}=\tau=0$. From formulas (\ref{2.1}) and (\ref{2.2}) we have
$u_{n\alpha}=0$ and $u_{\alpha\beta n}=0$. Then
\begin{align*}
W_{\alpha\beta n}  &  =u_{\alpha\beta n}+\frac{1-t}{n-2}\Sigma_{k}%
u_{kkn}g_{\alpha\beta}+au_{\alpha n}u_{\beta}+au_{\alpha}u_{\beta n}%
+a_{n}u_{\alpha}u_{\beta}\\
&  +2b\Sigma_{k}u_{kn}u_{k}g_{\alpha\beta}+b_{n}\left(  \Sigma_{k}u_{k}%
^{2}\right)  g_{\alpha\beta}+S_{\alpha\beta n}\\
&  =\frac{1-t}{n-2}u_{nnn}g_{\alpha\beta}+a_{n}u_{\alpha}u_{\beta}%
+b_{n}\left(  \Sigma_{\gamma}u_{\gamma}^{2}\right)  g_{\alpha\beta}%
+S_{\alpha\beta n}.
\end{align*}

For $V_{ij}$ we can get the equalities in the same way. \begin{flushright}
$\Box$
\end{flushright}

\noindent\textbf{Lemma 5.} \textit{Let }$(M^{n},g)$\textit{\ be a compact
Riemannian manifold with boundary and dimensional }$n\geq3$\textit{. Assume
that the boundary }$\partial M$\textit{\ is totally geodesic. Then at any
boundary point }$P\in\partial M$\textit{, there exists a conformal metric
}$\overline{g}=e^{-2\overline{u}}g$\textit{\ such that (i) }$\overline{u}%
_{n}=0$\textit{\ on }$\partial M$\textit{\ and the boundary }$\partial
M$\textit{\ is still totally geodesic,} \textit{(ii) }$\overline{R}_{ij}%
(P)=0$\textit{\ for }$1\leq i,j\leq n$\textit{, (iii) }$\overline{R}%
_{nn,n}(P)=0,\overline{R}_{\alpha n,\beta}(P)=0,1\leq\alpha,\beta\leq
n-1$\textit{, and (iv) }$\overline{R}_{\alpha\beta,n}(P)=0,1\leq\alpha
,\beta\leq n-1$\textit{.} \newline\textit{Proof.}\textbf{ } As the proof of
Lemma 3.3 in \cite{E2}, consider the first eigenvalue $\lambda_{1}\left(
L\right)  $ of the conformal Laplacian with respect to the boundary condition
\begin{equation}
\left\{
\begin{array}
[c]{ll}%
L\varphi+\lambda_{1}\left(  L\right)  \varphi=0 & \text{ \ \ \ \ on \ \ }M,\\
\frac{\partial}{\partial x^{n}}\varphi+\frac{n-2}{2}h\varphi=0 & \text{ \ on
\ \ }\partial M,
\end{array}
\right.  \tag{2.5}\label{2.5}%
\end{equation}
where $L=\Delta-\frac{n-2}{4\left(  n-1\right)  }R$, $R$ is the scalar
curvature, $h$ is the mean curvature of the boundary and $\frac{\partial
}{\partial x^{n}}$ is the inward norm derivative with respect to the metric
$g$. Since $\partial M$ is totally geodesic, $h=0$. Let $\varphi_{1}$ be the
first eigenfunction for the conformal Laplacian with respect to the boundary
condition (\ref{2.5}), then $\varphi_{1}>0$. Set $g_{1}=\varphi_{1}^{\frac
{4}{n-2}}g$. The transformation law of the second fundamental form
\[
\widetilde{L}_{\alpha\beta}=e^{f}L_{\alpha\beta}-\frac{\partial}{\partial
x^{n}}\left(  e^{f}\right)  g_{\alpha\beta}%
\]
with respect to the conformal change $\widetilde{g}=e^{2f}g$ implies that
$\partial M$ is totally geodesic. Recall $\left(  x^{1},...,x^{n-1}%
,x^{n}\right)  $ is Fermi coordinates around $P\in\partial M$. By Theorem 5.2
in \cite{LP}, there exists a homogeneous polynomial of degree 3, $k_{1}\left(
x\right)  $ such that the metric $g_{2}=e^{2k_{1}\left(  x\right)  }g_{1}$
satisfies $\frac{\partial k_{1}\left(  x\right)  }{\partial x^{n}}|_{\partial
M}=0$, $\partial M$ being totally geodesic, and
\[
R_{ij,k}\left(  P\right)  +R_{jk,i}\left(  P\right)  +R_{ki,j}\left(
P\right)  =0\text{ \ \ \ for all \ }1\leq i,j,k\leq n.
\]
We then get $R_{nn,n}\left(  P\right)  =0$ and $R_{\alpha\beta,n}\left(
P\right)  +R_{\beta n,\alpha}\left(  P\right)  +R_{n\alpha,\beta}\left(
P\right)  =0$ for $1\leq\alpha,\beta\leq n-1$. By the Codazzi equation for
$1\leq\alpha,\beta,\gamma\leq n-1$,
\begin{equation}
R_{\alpha\beta\gamma n}=L_{\alpha\gamma,\beta}-L_{\beta\gamma,\alpha}.
\tag{2.6}\label{2.6}%
\end{equation}
Differentiating (\ref{2.6}), we get for $1\leq\alpha,\beta,\gamma,\delta\leq
n-1$
\[
R_{\alpha\beta\gamma n,\delta}=L_{\alpha\gamma,\beta\delta}-L_{\beta
\gamma,\alpha\delta}.
\]
Since $\partial M$ is totally geodesic, after contracting with the metric, we
obtain for $1\leq\alpha,\beta\leq n-1$%
\[
R_{\alpha n,\beta}=0\text{ \ \ on \ }\partial M.
\]
Hence $R_{\alpha\beta,n}\left(  P\right)  =0$ for $1\leq\alpha,\beta\leq n-1$.
Let $\overline{g}=g_{2}=e^{2k_{1}\left(  x\right)  }\varphi_{1}^{\frac{4}%
{n-2}}g=e^{-2\overline{u}}g$, the metric $\overline{g}$ satisfies all the
properties we needed. \begin{flushright}
$\Box$
\end{flushright}

\section{Applications}

We denote $[g]=\{\hat{g}\mid\hat{g}=e^{-2u}g\}$ and $[g]_{k}=\{\hat{g}\mid
\hat{g}\in\lbrack g]\cap\Gamma_{k}^{+}\}.$ We call $g$ is $k$-admissible if
and only if $[g]_{k}\neq\emptyset$. Now the first Yamabe constant on
Riemannian manifold $\left(  M^{n},g\right)  $ with nonempty boundary
$\partial M$ can be defined as (\cite{E1})

\[
\mathcal{Y}_{1}[g]=\inf_{u\in C^{1}(M), u\neq0, \int_{M} u^{\frac{2n}{(n-2)}%
}=1} \left(  \int_{M}\left( |\nabla u|^{2}+\frac{n-2}{4(n-1)}R_{g}%
u^{2}\right)  +\frac{n-2}{2}\oint_{\partial M}h_{g}u^{2} \right)  .
\]
We may define the boundary curvature $\mathcal{B}^{k}$ for the manifold with
umbilic boundary and higher order Yamabe constants $\mathcal{Y}_{k}[g]$ for
$2\leq k<n/2$ as follows (these concepts were defined in \cite{Cn3}, the
similar higher order Yamabe constants for the manifolds without boundary have
appeared in \cite{GLW, S}, which are different with the "well-known"
definitions, e.g., see \cite{STW}):
\[
\mathcal{B}^{k}=\sum_{i=0}^{k-1}C(n,k,i)\sigma_{i}\left( \lambda( g^{-1}A^{T})\right)
\tau^{2k-2i-1},
\]
and
\[
\mathcal{Y}_{k}[g]=\left\{
\begin{array}
[c]{ll}%
inf_{\hat{g}\in\lbrack g]_{k-1};vol(\hat{g})=1}\mathcal{F}_{k} & \quad\text{if
\ }[g]_{k-1}\neq\emptyset\\
-\infty & \quad\text{if \ }[g]_{k-1}=\emptyset
\end{array}
\right.
\]
where $C(n,k,i)=\frac{\left(  n-i-1\right)  !}{\left(  n-k\right)  !\left(
2k-2i-1\right)  !!}$, $A^{T}=\left[  A_{\alpha\beta}\right]  $ is the
tangential part of the Schouten tensor, $\tau$ is a function satisfying
$L_{\alpha\beta}=\tau g_{\alpha\beta}$, and
\[
\mathcal{F}_{k}(\hat{g})=\int_{M}\sigma_{k}(\lambda(\hat{g}^{-1}A_{\hat{g}})
)+\oint_{\partial M}\mathcal{B}_{\hat{g}}^{k}.
\]
If $\partial M$ is totally geodesic with respect to $g$, then $\mathcal{B}_{g}^{k}=0$.
By Theorem 1 we can get the following Theorems 5 and 6 which can be viewed as
a generalization of the corresponding theorems in \cite{Cn3}.\\[0.2cm]
\textbf{Theorem 5. }\textit{Let }$(M,\partial M,g)$\textit{\ be a compact
manifold of dimension }$n\geq3$\textit{\ with boundary, }$\partial
M$\textit{\ is totally geodesic. Suppose that }$g\in [g]_{k-1}, 2\leq k<n/2$\textit{\ and
}$\mathcal{Y}_{1},\ \mathcal{Y}_{k}>0$\textit{. Then there exists a metric
}$\hat{g}\in\lbrack g]$\textit{\ such that }$A_{\hat{g}}\in\Gamma_{k}%
$\textit{\ and }$\partial M$ \textit{\ is totally geodesic under }$\hat{g}.$
\newline\textit{Proof.}\textbf{
}Following proof is mainly from \cite{GV1}. Comparing with \cite{S}, we may
prove the theorem by continuity method. Consider a family of equations
involving a parameter $t$,
\begin{equation}
\left\{
\begin{array}
[c]{ll}%
\sigma_k^{1/k}\left(\lambda(g^{-1}A_{\widehat{g}}^{t})\right)=f(x)e^{2u_{t}} & \text{ \ in
\ }M\text{\ }\\
\frac{\partial u_{t}}{\partial\nu}=0 & \text{ \ on \ }\partial M
\end{array}
\right.  \tag{3.1}\label{3.1}%
\end{equation}
where $\widehat{g}=e^{-2u_{t}}g$, $f\left(  x\right)  >0$ and $t\leq1$. Since
$g\in [g]_{k-1}$, the scalar curvature $R_{g}>0$. Then there
exists $a>-\infty$ so that $A_{g}^{a}$ is positive definite. For $t\in\left[
a,1\right]  $, we consider the deformation
\begin{equation}
\left\{
\begin{array}
[c]{ll}%
\sigma_k^{1/k}\left(\lambda\left(g^{-1}A_{u_{t}}^{t}\right)\right)=f(x)e^{2u_{t}} & \text{ \ in \ }M\text{\ }\\
\frac{\partial u_{t}}{\partial\nu}=0 & \text{ \ on \ }\partial M
\end{array}
\right.  \tag{3.2}\label{3.2}%
\end{equation}
where $A_{u_{t}}^{t}=A_{\widehat{g}}^{t}$ with $\widehat{g}=e^{-2u_{t}}g$,
$f\left(  x\right)  =\sigma_k^{1/k}\left( \lambda( g^{-1}A_{u_{a}}^{a})\right)  >0$ and
$u_{a}\equiv0$ is a solution of (\ref{3.2}) for $t=a$. Let
\[
I=\left\{  t\in\left[  a,t_{0}\right]  |%
\begin{array}
[c]{c}%
\exists\text{ a solution }u_{t}\in C^{2,\alpha}\left(  M\right)  \text{ of
(\ref{3.2}) with }{A_{\hat{g}_{t}}\in\Gamma_{k}}\text{ }\\
\text{and }\partial M\text{ being totally geodesic under }\hat{g}_{t}%
\end{array}
\right\}  .
\]
It is easy to prove that the linearized operator $L_{t}:C^{2,\alpha}\left(
M\right)  \cap\left\{  \frac{\partial u}{\partial\nu}|_{\partial M}=0\right\}
\rightarrow C^{\alpha}\left(  M\right)  $ is invertible. This together with
the implicit theorem imply that the set $I$ is open.

Theorem 1 implies the $C^{1}$ and $C^{2}$ estimates of the solution to
(\ref{3.2}) which depend only on the upper bound of $u$. Since $A^{t}%
=A^{1}+\frac{1-t}{n-2}\sigma_{1}\left(  A^{1}\right)  g$, at the maximal point
$x_{0}$ of $u_{t}$, we have $\left\vert \nabla u_{t}\right\vert =0$ and
$\nabla^{2}u_{t}\left(  x_{0}\right)  $ is negative semi-definite, no matter
$x_{0}$ being interior or boundary point. Hence,
\[
f(x_{0})^{k}e^{2ku\left(  x_{0}\right)  }=\sigma_{k}\left( \lambda( g^{-1}A_{u_{t}%
}^{t})\right)  \leq\sigma_{k}\left( \lambda\left( g^{-1}\left(  A+\frac{1-t}{n-2}\sigma
_{1}\left(  A\right)  g\right) \right) \right)  \leq C,
\]
where we use $\sigma_{1}\left(  A\right)  >0$ and $a\leq t\leq1$. We then get
the upper bound. By the gradient estimate and the assumption $\Y_1>0, \Y_k>0$,
we may easily get the lower bound
of $u$. Therefore we conclude that $I=\left[  a,1\right]$. We thus finish
the proof. \begin{flushright} $\Box$
\end{flushright}

If {$(M,g)$ is a locally conformally flat compact manifold of dimension
$n\geq3$ with umbilic boundary. Then by \cite{E1}, we may assume that the
background metric }$g$ is a Yamabe metric with its constant scalar curvature
$R>0$ and the boundary is totally geodesic. Then using the same argument of
Theorem 5, we may prove that {there exists a metric $\hat{g}\in\lbrack g]$
such that $A_{\hat{g}}\in\Gamma_{k}$ and $\partial M$ is totally geodesic under $\hat{g}$. By
\cite{JLL}, we can get the following existence result.\\[0.2cm]
\textbf{Theorem 6. }\textit{\ Let }$(M,\partial M,g)$\textit{\ be a locally
conformally flat compact manifold of dimension }$n\geq3$\textit{\ with umbilic
boundary. Suppose that }$2\leq k<n/2$\textit{\ and }$\Y_{1},\ \Y_{k}>0$\textit{.
Then there exists a metric }$\hat{g}\in\lbrack g]$\textit{\ such that }%
$\sigma_{k}(\lambda(A_{\hat{g}}))=1$\textit{\ and }$\partial M$\textit{\ is totally geodesic under }$\hat{g}.$

\section{Proof of Theorems 1 and 2}

\noindent In this section we first prove Theorem 2 for some general functions
$a\left(  x\right)  $ and $b\left(  x\right)  $, and a general $2$-symmetric
tensor $S$. After we establish the interior a priori estimates (Theorem 2), we
study the boundary estimates for special functions $a\left(  x\right)  $,
$b\left(  x\right)  $ and a special $2$-symmetric tensor $S$.

\noindent\textit{Proof of Theorem 2.}

\noindent\textit{(1) Case (a)}.

Let $K=\triangle u+a|\nabla u|^{2}$. Note that $\Gamma\subset\Gamma_{1}^{+}$ ,
we can immediately get
\[
0\leq\ tr(W)=(1+n\frac{1-t}{n-2})\triangle u+(a+nb)|\nabla u|^{2}%
+trS\leq(1+n\frac{1-t}{n-2})K-n\delta_{1}|\nabla u|^{2}+C.
\]
Then%
\[
(1+n\frac{1-t}{n-2})K\geq n\delta_{1}|\nabla u|^{2}-C>-C.
\]
Hence, $K$ has lower bound. We also have
\begin{equation}
|\nabla u|^{2}\leq\frac{(1+n\frac{1-t}{n-2})K+C}{n\delta_{1}}. \tag{4.1}%
\label{4.1}%
\end{equation}
Without loss of generality, we may assume $K>0$. Otherwise, $K\le0$. By the
above inequality (\ref{4.1}), we know that $|\nabla u|^{2}\leq C$. Then we
have the $C^{1}$ estimates. Furthermore, we have $|\triangle u|\leq C$. From
the condition $\Gamma\subset\Gamma_{2}^{+}$, we know that $\left(  trW\right)
^{2}-\left\vert W\right\vert ^{2}=2\sigma_{2}\left(  W\right)  >0$. Therefore
$\left\vert W\right\vert \leq Ctr(W)\leq C$ which implies $\left\vert
\nabla^{2}u\right\vert \leq C.$ We then get $C^{2}$ estimates.

Now by the assumption and (\ref{4.1}), we have
\begin{equation}
|\nabla u|^{2}\leq C(K+1), \tag{4.2}\label{4.2}%
\end{equation}
where $C$ depends on $||a||_{\infty}$ and $||b||_{\infty}$. By (\ref{4.2}), we
can obtain
\[
\triangle u=K-a|\nabla u|^{2}\leq K+||a||_{\infty}|\nabla u|^{2}\leq C(K+1).
\]
By the condition $\Gamma\subset\Gamma_{2}^{+}$ again, we know that
$|W_{ij}|\leq Ctr(W)$ which implies
\begin{equation}
|\nabla^{2}u|\leq C(K+1), \tag{4.3}\label{4.3}%
\end{equation}
where $C$ depends only on $||a||_{\infty}$ and $||b||_{\infty}$ as well.
(\ref{4.2}) and (\ref{4.3}) are the fundamental inequalities which we will use
over and over again.

In order to prove that $K$ is bounded, similar as \cite{Cn2}, we consider an
auxiliary function $H=\eta K$ in a neighborhood $B_{r}$, where $0\leq\eta
\leq1$ is a cutoff function depending only on $r$ such that $\eta=1$ in
$B_{\frac{r}{2}}$ and $\eta=0$ outside $B_{r}$, $|\nabla\eta|\leq\frac
{C\eta^{1/2}}{r},|\nabla^{2}\eta|\leq\frac{C}{r^{2}}$.

We begin to derive the interior $C^{1}$ and $C^{2}$ estimates.

At the maximum point of $H$, $x_{0}$, after choosing normal coordinates, we
have
\[
0=H_{i}=\eta_{i}K+\eta K_{i}%
\]
That is
\[
K_{i}=-\frac{\eta_{i}}{\eta}K.
\]
We also have
\[
0\geq H_{ij}=\eta_{ij}K+\eta K_{ij}+\eta_{i}K_{j}+\eta_{j}K_{i}.
\]
Note that $|\nabla\eta|\leq\frac{C\eta^{1/2}}{r},|\nabla^{2}\eta|\leq\frac
{C}{r^{2}}$, we have
\[
0\geq H_{ij}=\eta K_{ij}+\Lambda_{ij}K,
\]
where
\[
\Lambda_{ij}=\eta_{ij}-\frac{2\eta_{i}\eta_{j}}{\eta}\geq-C\delta_{ij}%
\]
and $C$ depends only on $r$.

Denote
\[
P^{ij}=F^{ij}+\frac{1-t}{n-2}\left(  \Sigma_{l}F^{ll}\right)  \delta^{ij}.
\]
Since $t\le 1$, $P^{ij}$ is still elliptic. By use of Ricci identities, we have
\[|u_{ijkk}-u_{kkij}|\le C|\nabla^2 u| \le C(K+1),\]
and
\[|u_{ijk}-u_{kij}|\le C|\nabla u| \le C(K^{1/2}+1).\]
We then have
\begin{align}
0  &  \geq\eta P^{ij}H_{ij}\nonumber\\
&  =\eta^{2}P^{ij}K_{ij}+\eta\Lambda_{ij}P^{ij}K\nonumber\\
&  \geq\eta^{2}P^{ij}\sum_{k}[u_{ijkk}+2a(u_{ki}u_{kj}+u_{ijk}u_{k}%
)+a_{ij}u_{k}^{2}+4a_{i}u_{kj}u_{k}]\nonumber\\
&  -C\left(  \sum_{i}F^{ii}\right)  (1+K)\nonumber\\
&  \geq\eta^{2}P^{ij}\sum_{k}[u_{ijkk}+2a(u_{ki}u_{kj}+u_{ijk}u_{k}%
)]-C(\sum_{i}F^{ii})(1+\eta^{3/2}K^{3/2}). \tag{4.4}\label{4.4}%
\end{align}
Now we estimate the terms $\sum_{k}P^{ij}u_{ijkk}$ and $\sum_{k}P^{ij}%
(u_{ki}u_{kj}+u_{ijk}u_{k})$ respectively.%
\begin{align*}
\sum_{k}P^{ij}u_{ijkk}  &  =F^{ij}\sum_{k}[W_{ijkk}-2a(u_{ik}u_{jk}%
+u_{ikk}u_{j})-2b(u_{lkk}u_{l}+u_{lk}u_{lk})g_{ij}\allowdisplaybreaks\\
&  -a_{kk}u_{i}u_{j}-4a_{k}u_{ik}u_{j}-b_{kk}u_{l}^{2}g_{ij}-4b_{k}u_{lk}%
u_{l}g_{ij}-S_{ijkk}]\allowdisplaybreaks\\
&  \geq\sum_{k}f_{kk}+F^{ij}\sum_{k}[-2a(u_{ik}u_{jk}+u_{kki}u_{j}%
)-2b(u_{kkl}u_{l}+u_{lk}u_{lk})g_{ij}]\allowdisplaybreaks\\
&  -C\sum F^{ii}(1+K^{3/2})\allowdisplaybreaks\\
&  \geq\sum_{k}f_{kk}+F^{ij}\sum_{k}[-2au_{ik}u_{jk}-2au_{j}(-2au_{ki}%
u_{k}-a_{i}u_{k}^{2}-\frac{\eta_{i}}{\eta}K)\allowdisplaybreaks\\
&  -2bu_{lk}^{2}g_{ij}-2bu_{l}(-2au_{kl}u_{k}-a_{l}u_{k}^{2}-\frac{\eta_{l}%
}{\eta}K)g_{ij}]\allowdisplaybreaks\\
&  -C\sum F^{ii}(1+K^{3/2}).
\end{align*}
We then get
\begin{align}
\sum_{k}P^{ij}u_{ijkk}  &  \geq &  &  \sum_{k}f_{kk}+F^{ij}\sum_{k}%
[-2b\Sigma_{l}\left(  u_{lk}\right)  ^{2}g_{ij}-2au_{ik}u_{jk}+4a^{2}%
u_{j}u_{ki}u_{k}\nonumber\\
&  &  &  +4ab\Sigma_{l}\left(  u_{l}u_{kl}u_{k}\right)  g_{ij}]-C\eta
^{-1/2}\sum F^{ii}(1+K^{3/2}). \tag{4.5}\label{4.5}%
\end{align}
We also have
\begin{align}
2aP^{ij}(u_{ijk}u_{k}+u_{ki}u_{kj})  &  = &  &  2aF^{ij}(W_{ijk}u_{k}%
-2au_{ik}u_{j}u_{k}-2bu_{lk}u_{l}g_{ij}u_{k}+u_{ik}u_{jk}\nonumber\\
&  &  &  +\frac{1-t}{n-2}\sum_{l,k}u_{lk}^{2}\delta_{ij}-a_{k}u_{i}u_{j}%
u_{k}-b_{k}u_{l}^{2}g_{ij}u_{k}-S_{ijk}u_{k})\nonumber\\
&  \geq &  &  2au_{k}f_{k}+F^{ij}(-4a^{2}u_{ik}u_{j}u_{k}-4abu_{lk}u_{l}%
u_{k}g_{ij}+2au_{ik}u_{jk})\nonumber\\
&  &  &  +\frac{1-t}{n-2}2a\sum_{i}F^{ii}\sum_{l,k}u_{lk}^{2}-C\sum_{i}%
F^{ii}(1+K^{3/2}). \tag{4.6}\label{4.6}%
\end{align}
From (\ref{4.4}), (\ref{4.5}) and (\ref{4.6}) we therefore have%
\begin{equation}
0\geq\left(  \sum F^{ii}\right)  (2(\frac{1-t}{n-2}a-b)\eta^{2}|\nabla
^{2}u|^{2}-C\eta^{3/2}K^{3/2}-C). \tag{4.7}\label{4.7}%
\end{equation}

Let $A$ be a number such that $A>
\sqrt{\frac{2}{\delta_2}}(\frac{1-t}{n-2})$. First, we assume  $|\nabla u|^{2}({x_{0}
})<A|\triangle u|({x_{0}})$. By $|u_{ij}|\leq C(K+1)$, we know that at the
point $x_{0}$, $|u_{ij}|\leq C(|\triangle u|+1)$. Thus (\ref{4.7}) becomes
\[
0\geq \sum F^{ii}(\frac{2}{n}(\frac{1-t}{n-2}a-b)\eta
^{2}|\triangle u|^{2}-C\eta^{3/2}|\triangle u|^{3/2}-C).
\]
Hence%
\[
|\triangle u|({x_{0}})\leq C
\]
and
\[
K\leq C.
\]
Next we consider the case $|\nabla u|^{2}({x_{0}})\geq
A|\triangle u|({x_{0}})$\textbf{.} From $|u_{ij}|\leq C(K+1)$, we know that at
the point $x_{0}$, $|u_{ij}|\leq C(|\nabla u|^{2}+1)$. Thus (\ref{4.7})
becomes
\begin{equation}
0\geq\sum F^{ii}(2(\frac{1-t}{n-2}a-b)\eta^{2}u_{il}^{2}-C\eta^{3/2}|\nabla
u|^{3}-C). \tag{4.8}\label{4.8}%
\end{equation}
We may assume that $W_{ij}$ is diagonal at the point $x_{0}$,
\[
W_{ii}=u_{ii}+\frac{1-t}{n-2}\triangle u+au_{i}^{2}+b\Sigma_{k}u_{k}%
^{2}+S_{ii},
\]
and
\[
0=W_{ij}=u_{ij}+au_{i}u_{j}+S_{ij},\ (i\neq j).
\]
Since
\[
F^{ii}(u_{il}+\frac{1-t}{n-2}\triangle ug_{il}+S_{il})^{2}\leq2F^{ii}%
[u_{il}^{2}+(\frac{1-t}{n-2}\triangle ug_{il}+S_{il})^{2}],
\]
we obtain%
\begin{align*}
2\Sigma_{l,i}F^{ii}u_{li}u_{li}  &  \geq\Sigma_{l,i}F^{ii}[(u_{il}+\frac
{1-t}{n-2}\triangle ug_{il}+S_{il})^{2}-2(\frac{1-t}{n-2}\triangle
ug_{il}+S_{il})^{2}]\\
&  \geq\Sigma_{i}F^{ii}[\sum_{j\neq i}(-au_{i}u_{j})^{2}+(W_{ii}-au_{i}%
^{2}-b|\nabla u|^{2})^{2}]\\
&  -2\left(  \frac{1-t}{n-2}\right)  ^{2}\frac{1}{A^{2}}\Sigma_{i}%
F^{ii}|\nabla u|^{4}-C\Sigma_{i}F^{ii}\\
&  \geq\Sigma_{i}F^{ii}[(au_{i})^{2}|\nabla u|^{2}+W_{ii}^{2}-2W_{ii}%
(au_{i}^{2}+b|\nabla u|^{2})+b^{2}|\nabla u|^{4}\\
&  +2abu_{i}^{2}|\nabla u|^{2}]-2\left(  \frac{1-t}{n-2}\right)  ^{2}\frac
{1}{A^{2}}\Sigma_{i}F^{ii}|\nabla u|^{4}-C\Sigma_{i}F^{ii}.
\end{align*}
Since
\[
2aF^{ii}W_{ii}u_i^2\le F^{ii}(W_{ii}^2+a^2u_i^4)\le F^{ii}W_{ii}^2+a^2F^{ii}u_i^2|\nabla u|^2,
\]
we have
\begin{align*}
2\Sigma_{l,i}F^{ii}u_{li}u_{li}  &  \geq \Sigma_{i}F^{ii}[(2ab+b^{2})u_{i}^{2}+\sum_{j\neq i}b^{2}u_{j}
^{2}]|\nabla u|^{2}-2||b||_{\infty}f|\nabla u|^{2}\\
&  -2\left(  \frac{1-t}{n-2}\right)  ^{2}\frac{1}{A^{2}}\Sigma_{i}
F^{ii}|\nabla u|^{4}-C\Sigma_{i}F^{ii}.
\end{align*}

By the assumption of the theorem case (a), $\min\left\{  2ab+b^{2}, b^{2}\right\}  \geq\delta_{2}>0,$ and Lemma 1, we then have
\[
2\Sigma_{l,i}F^{ii}u_{li}u_{li}\geq\left(  \delta_{2}-2\left(  \frac{1-t}%
{n-2}\right)  ^{2}\frac{1}{A^{2}}\right)  \Sigma_{i}F^{ii}|\nabla
u|^{4}-C\Sigma_{i}F^{ii}|\nabla u|^2-C\Sigma_{i}F^{ii}
\]
Therefore, by (\ref{4.8}) we have
\[
0\geq\left(  \Sigma_{i}F^{ii}\right)  [\delta_{1}\left(  \delta_{2}-2\left(
\frac{1-t}{n-2}\right)  ^{2}\frac{1}{A^{2}}\right)  \eta^{2}|\nabla
u|^{4}-C\eta^{3/2}|\nabla u|^{3}-C\eta|\nabla u|^{2}-C].
\]
Since $A>0$ large enough, we have
\[
|\nabla u|^{2}({x_{0}})\leq C,
\]
therefore
\[
K\leq C.
\]

\noindent\textit{Proof of Remark.} We may estimate the term $2aF^{ii}W_{ii}u_i^2$ as follows. Since $W\in \Gamma_2^+$, we have
\[
W_{ii}\le tr W=(1+n\frac{1-t}{n-2})\triangle u+(a+nb)|\nabla u|^2+tr S
\]
for each $i$. By use of the condition $|\nabla u|^2(x_0)\ge A\triangle u(x_0)$ for some suitable large number $A$, we have
\[
W_{ii}\le [(1+n\frac{1-t}{n-2})\frac{1}{A}+(a+nb)]|\nabla u|^2+tr S.
\]
Now
\[
2aF^{ii}W_{ii}u_i^2\le 2||a||_{\infty}F^{ii}u_i^2\{[(1+n\frac{1-t}{n-2})\frac{1}{A}+(a+nb)]|\nabla u|^2+tr S\}.
\]
Then
\begin{align*}
2\Sigma_{l,i}F^{ii}u_{li}u_{li}  &  \geq \sum_{i}F^{ii}[ \Big( b^2 -a^{2}+2ab-2n||a||_{\infty}b
\Big)u_{i}^{2}+\sum_{j\neq i}b^{2}u_{j} ^{2}]|\nabla u|^{2}
  \\
  &  -2||b||_{\infty}f|\nabla u|^{2}
 -2\left(  \frac{1-t}{n-2}\right)  ^{2}\frac{1}{A^{2}}\sum_{i}
F^{ii}|\nabla u|^{4} \\
& -2||a||_{\infty}(1+n\frac{1-t}{n-2})\frac{1}{A}\sum_{i}F^{ii} u_{i}^{2}
|\nabla u|^2-C\sum_{i}F^{ii}
\end{align*}
Using the condition $\min \{b^2-a^2+2ab-2n||a||_{\infty}b, b^2\}\ge \delta_2>0$ in Remark, we may get
\begin{align*}
2\Sigma_{l,i}F^{ii}u_{li}u_{li}\geq & \left(  \delta_{2}-2\left(  \frac{1-t}
{n-2}\right)  ^{2}\frac{1}{A^{2}}-2||a||_{\infty}(1+n\frac{1-t}{n-2})\frac{1}{A}\right)  \Sigma_{i}F^{ii}|\nabla
u|^{4}\\
&-C\Sigma_{i}F^{ii}|\nabla u|^2-C\Sigma_{i}F^{ii}.
\end{align*}
Now substituting this inequality to (\ref{4.8}) we may get the desire estimate.

\noindent\textit{(2) Case (b).}

Since $a(x)+nb(x)\le -\delta_{3},$ by the condition $\Gamma\subset\Gamma_{1}^{+}%
$, we have
\begin{align*}
0  &  \leq\ tr(W)=(1+n\frac{1-t}{n-2})\triangle u+a|\nabla u|^{2}+nb|\nabla
u|^{2}+tr S\\
&  \leq(1+n\frac{1-t}{n-2})\triangle u-\delta_{3}|\nabla u|^{2}+C.
\end{align*}
Then
\begin{equation}
|\nabla u|^{2}\leq C(\triangle u+1). \tag{4.9}\label{4.9}%
\end{equation}

The proof is similar as the argument in case (a). We take the same auxiliary
function $H=\eta(\triangle u+a|\nabla u|^{2})\triangleq\eta K$, where
$\eta\left(  r\right)  $ is a cutoff function as in case (a). Without loss of
generality, we may assume
\[
K=\triangle u+a|\nabla u|^{2}>>1.
\]
Since $a\left(  x\right)  \geq0$, by (\ref{4.9}), we have
\begin{equation}
\Delta u\leq C\left(  K+1\right)  \tag{4.10}\label{4.10}%
\end{equation}
and
\begin{equation}
|\nabla u|^{2}\leq C(K+1). \tag{4.11}\label{4.11}%
\end{equation}

Suppose that the maximum point of $H$ achieves at $x_{0}$, an interior point.
Then at this point, we need to note that $|\nabla u|^{2},\triangle u$ and $K$
all can be controlled by $C\left(  \left\vert \nabla^{2}u\right\vert
+1\right)  $. By the same computation as in case (a), (\ref{4.4}),
(\ref{4.5}), (\ref{4.6}) and (\ref{4.7}) become
\begin{equation}
0\geq\eta^{2}P^{ij}\sum_{k}[u_{ijkk}+2a(u_{ki}u_{kj}+u_{ijk}u_{k})]-C(\sum
_{i}F^{ii})(1+\left\vert \nabla^{2}u\right\vert ^{3/2}), \tag{4.12}%
\label{4.12}%
\end{equation}%
\begin{align}
\sum_{k}P^{ij}u_{ijkk}  &  \geq &  &  \sum_{k}f_{kk}+F^{ij}\sum_{k}%
[-2b\Sigma_{l}\left(  u_{lk}\right)  ^{2}g_{ij}-2au_{ik}u_{jk}+4a^{2}%
u_{j}u_{ki}u_{k}\nonumber\\
&  &  &  +4ab\Sigma_{l}\left(  u_{l}u_{kl}u_{k}\right)  g_{ij}]-C\eta
^{-1/2}\sum F^{ii}(1+\left\vert \nabla^{2}u\right\vert ^{3/2}), \tag{4.13}%
\label{4.13}%
\end{align}%
\begin{align}
2aP^{ij}(u_{ijk}u_{k}+u_{ki}u_{kj})  &  \geq &  &  2au_{k}f_{k}+F^{ij}%
(-4a^{2}u_{ik}u_{j}u_{k}-4abu_{lk}u_{l}u_{k}g_{ij}+2au_{ik}u_{jk})\nonumber\\
&  &  &  +\frac{1-t}{n-2}2a\sum_{i}F^{ii}\sum_{l,k}u_{lk}^{2}-C\sum_{i}%
F^{ii}(1+\left\vert \nabla^{2}u\right\vert ^{3/2}), \tag{4.14}\label{4.14}%
\end{align}
and
\begin{equation}
0\geq\left(  \sum F^{ii}\right)  (2(\frac{1-t}{n-2}a-b)\left(  \eta|\nabla
^{2}u|\right)  ^{2}-C\left(  \eta\left\vert \nabla^{2}u\right\vert \right)
^{3/2}-C\eta\left\vert \nabla^{2}u\right\vert -C) \tag{4.15}\label{4.15}%
\end{equation}
respectively. (\ref{4.15}) gives $\eta|\nabla^{2}u|\left(  x_{0}\right)  \leq
C$ and hence the bounds of $K,|\nabla^{2}u|$ and $\left\vert \nabla
u\right\vert $. \noindent\begin{flushright} $\Box$
\end{flushright}

\noindent\textit{Proof of Theorem 1.}

Note in this theorem $a,b$ are constants and $S=A$ is the Schouten tensor.
Similar with the case (b) in Theorem 2. Let $K=\triangle u+a|\nabla u|^{2}$.
We have (\ref{4.9}) , (\ref{4.10}) and (\ref{4.11}). Consider $\overline
{H}=\eta Ke^{px_{n}}$ in a neighborhood $\overline{B}^{+}_{r}$, where
$0\leq\eta\leq1$ is a cutoff function depending only on $r$ such that $\eta=1$
in $\overline{B}^{+}_{\frac{r}{2}}$ and $\eta=0$ outside $\overline{B}^{+}%
_{r}$, and $p$ is a large positive constant.

\textit{Step 1.} We first prove the maximum point of $\overline{H}$ must be in
the interior of $M$. Assume $\overline{H}$ arrives at its maximum point
$x_{0}$ on the boundary. By a direct calculation of $\overline{H}_{n}$, we can
show that $\overline{H}_{n}|_{x_{0}}>0$, which violates the assumption.

Note that $\eta$ is a function depending only on $r$, thus at the boundary
point we have $\eta_{n}=0$. By use of (\ref{2.1}) and (\ref{2.2}), we get%
\begin{align}
\overline{H}_{n}|_{x_{0}}  &  =\eta e^{px_{n}}(K_{n}+pK)\nonumber\\
&  =\eta e^{px_{n}}(u_{nnn}+u_{\alpha\alpha n}+2au_{\alpha n}u_{\alpha
}+2au_{nn}u_{n}+a_{n}(u_{\gamma}u_{\gamma}+u_{n}u_{n})+pK)\nonumber\\
&  =\eta e^{px_{n}}(u_{nnn}+pK), \tag{4.16}\label{4.16}%
\end{align}
where the last equality follows from
\begin{align*}
&  u_{\alpha\alpha n}+2au_{\alpha n}u_{\alpha}+2au_{nn}u_{n}+a_{n}(u_{\gamma
}u_{\gamma}+u_{n}u_{n})\\
&  =a_{n}u_{\gamma}u_{\gamma}=0.
\end{align*}

Now we need the following Lemma 6:

\noindent\textbf{Lemma 6. }\ \textit{There exists a constant }%
$\mathit{\mathit{C}}$\textit{ depending only on }$\mathit{a,b,\Lambda,t}%
$\textit{ and $\varepsilon$, such that at $x_{0}$, $u_{nnn}\geq-C(K+1).$}

By Lemma 6 we know that if $p$ is large enough,
\begin{align*}
\overline{H}_{n}|_{x_{0}}  &  =\eta e^{px_{n}}(u_{nnn}+pK)\\
&  \geq\eta e^{px_{n}}((p-C)K-C)>0,
\end{align*}
which completes the first step of the proof.

\noindent\textit{Proof of Lemma 6}.\textbf{ }

By Lemma 5, we may choose a conformal metric $\bar{g}=e^{-2\bar{u}}g$ and
$\overline{u}_{n}|_{\partial M}=0$ at first. In this metric, $\partial M$ is
still totally geodesic and $\overline{A}_{\alpha\beta,n}(x_{0})=0$. We wish to
find a metric $\widetilde{g}=e^{-2v}\overline{g}$ such that $u=\overline{u}+v$
is a solution to (\ref{1.5}). Now
\begin{align*}
W_{i}^{l}  &  =g^{lj}\left(  \frac{1-t}{n-2}\Delta ug_{ij}+u_{ij}+au_{i}%
u_{j}+b\left\vert \nabla u\right\vert ^{2}g_{ij}+A_{ij}\right) \\
&  =g^{lj}\left(  \left(  \overline{u}+v\right)  _{ij}+a\left(  \overline
{u}+v\right)  _{i}\left(  \overline{u}+v\right)  _{j}+b\left\vert
\nabla\left(  \overline{u}+v\right)  \right\vert ^{2}g_{ij}+A_{ij}\right) \\
&  +\frac{1-t}{n-2}\Delta\left(  \overline{u}+v\right)  \delta_{i}^{l}\\
&  =g^{lj}\left(  \overline{A}_{ij}+\left(  a-1\right)  \overline{u}%
_{i}\overline{u}_{j}+\left(  b+\frac{1}{2}\right)  \left\vert \nabla
\overline{u}\right\vert ^{2}g_{ij}\right)  +\frac{1-t}{n-2}\Delta\left(
\overline{u}+v\right)  \delta_{i}^{l}\\
&  +g^{lj}\left(  v_{ij}+a\overline{u}_{i}v_{j}+a\overline{u}_{j}v_{i}%
+av_{i}v_{j}+b\left(  \overline{u}_{k}v_{l}+\overline{u}_{l}v_{k}\right)
g^{kl}g_{ij}+b\left\vert \nabla v\right\vert ^{2}g_{ij}\right) \\
&  =e^{-2\overline{u}}\overline{g}^{lj}\left(  \overline{A}_{ij}+\left(
a-1\right)  \overline{u}_{i}\overline{u}_{j}+\left(  b+\frac{1}{2}\right)
\left\vert \overline{\nabla}\overline{u}\right\vert _{\overline{g}}%
^{2}\overline{g}_{ij}\right) \\
&  +e^{-2\overline{u}}\overline{g}^{lj}(\overline{\nabla}_{ij}^{2}v+\left(
\overline{\Gamma}_{ij}^{k}\left(  \overline{g}\right)  -\Gamma_{ij}^{k}\left(
g\right)  \right)  v_{k}+a\overline{u}_{i}v_{j}+a\overline{u}_{j}v_{i}%
+av_{i}v_{j})\\
&  +e^{-2\overline{u}}\overline{g}^{lj}\left(  b\left(  \overline{u}_{k}%
v_{l}+\overline{u}_{l}v_{k}\right)  \overline{g}^{kl}\overline{g}%
_{ij}+b\left\vert \overline{\nabla}v\right\vert _{\overline{g}}^{2}%
\overline{g}_{ij}\right) \\
&  +e^{-2\overline{u}}\frac{1-t}{n-2}\left(  \overline{\Delta}\left(
\overline{u}+v\right)  +\overline{g}^{jk}\left(  \overline{\Gamma}_{ij}%
^{k}\left(  \overline{g}\right)  -\Gamma_{ij}^{k}\left(  g\right)  \right)
\left(  \overline{u}_{k}+v_{k}\right)  \right)  \delta_{i}^{l},
\end{align*}
where $\overline{A}_{ij}=\overline{u}_{ij}+\overline{u}_{i}\overline{u}%
_{j}-\frac{1}{2}\left\vert \nabla\overline{u}\right\vert ^{2}g_{ij}+A_{ij}.$
Then equation (\ref{1.5}) becomes
\begin{equation}
\left\{
\begin{array}
[c]{ll}%
F(\bar{g}^{-1}\overline{W})=e^{2\overline{u}}f(x,\overline{u}+v) & \text{
in\ }\overline{B}_{r}^{+},\\
\frac{\partial v}{\partial x^{n}}=0 & \text{ on }\Sigma_{r},
\end{array}
\right.  \tag{4.17}\label{4.17}%
\end{equation}
where
\begin{align*}
\overline{W}_{ij}  &  =\overline{A}_{ij}+\left(  a-1\right)  \overline{u}%
_{i}\overline{u}_{j}+\left(  b+\frac{1}{2}\right)  \left\vert \overline
{\nabla}\overline{u}\right\vert _{\overline{g}}^{2}\overline{g}_{ij}\\
&  +\overline{\nabla}_{ij}^{2}v+\left(  \overline{\Gamma}_{ij}^{k}\left(
\overline{g}\right)  -\Gamma_{ij}^{k}\left(  g\right)  \right)  v_{k}%
+a\overline{u}_{i}v_{j}+a\overline{u}_{j}v_{i}+av_{i}v_{j}\\
&  +b\left(  \overline{u}_{k}v_{l}+\overline{u}_{l}v_{k}\right)  \overline
{g}^{kl}\overline{g}_{ij}+b\left\vert \overline{\nabla}v\right\vert
_{\overline{g}}^{2}\overline{g}_{ij}\\
&  +\frac{1-t}{n-2}\left(  \overline{\Delta}\left(  \overline{u}+v\right)
+\overline{g}^{lk}\left(  \overline{\Gamma}_{lk}^{p}\left(  \overline
{g}\right)  -\Gamma_{lk}^{p}\left(  g\right)  \right)  \left(  \overline
{u}_{p}+v_{p}\right)  \right)  \overline{g}_{ij}.
\end{align*}
Since $u_{n}=\overline{u}_{n}=0$ on the boundary $\partial M$, $u_{n\alpha
}=\overline{u}_{n\alpha}=0$. By Lemma 3, we can show $u_{\alpha\beta
n}=\overline{u}_{a\beta n}=0$, therefore $v_{n}=v_{n\alpha}=v_{\alpha\beta
n}=0$ on $\partial M$. We then have
\[
\overline{W}_{\alpha n}(x_{0})=\overline{A}_{\alpha n}+v_{\alpha n}+\left(
\overline{\Gamma}_{\alpha n}^{\beta}\left(  \overline{g}\right)
-\Gamma_{\alpha n}^{\beta}\left(  g\right)  \right)  v_{\beta}=0.
\]
Applying an argument of Lemma 13 in \cite{Cn3}, we know $F^{\alpha n}%
(x_{0})=0$. Now by Lemma 5,
\begin{align*}
\overline{W}_{\alpha\beta n}\left(  x_{0}\right)   &  =\left(  \overline
{\Gamma}_{\alpha\beta}^{\delta}\left(  \overline{g}\right)  -\Gamma
_{\alpha\beta}^{\delta}\left(  g\right)  \right)  _{n}v_{\delta}|_{x_{0}}\\
&  +\frac{1-t}{n-2}\left(  \overline{u}_{nnn}+v_{nnn}+\left[  \overline
{g}^{\gamma\delta}\left(  \overline{\Gamma}_{\gamma\delta}^{p}\left(
\overline{g}\right)  -\Gamma_{\gamma\delta}^{p}\left(  g\right)  \right)
\left(  \overline{u}_{p}+v_{p}\right)  \right]  _{,n}\right)  \overline
{g}_{\alpha\beta}.
\end{align*}
Here we have used the fact $\overline{g}_{ij,n}=\overline{g}_{,n}^{ij}=0$. Use
Fermi coordinates, we have on $\partial M$
\[
\frac{\partial g_{\alpha\beta}}{\partial x^{n}}=\frac{\partial}{\partial
x^{n}}<\frac{\partial}{\partial x^{\alpha}},\frac{\partial}{\partial x^{\beta
}}>=<\nabla_{\frac{\partial}{\partial x^{\alpha}}}\frac{\partial}{\partial
x^{n}},\frac{\partial}{\partial x^{\beta}}>+<\nabla_{\frac{\partial}{\partial
x^{\beta}}}\frac{\partial}{\partial x^{n}},\frac{\partial}{\partial x^{\alpha
}}>=-2L_{\alpha\beta},
\]
and
\begin{align*}
\frac{\partial}{\partial x^{n}}\Gamma_{\alpha\beta}^{\delta}\left(  g\right)
&  =\frac{1}{2}g^{\delta\gamma}\left(  \frac{\partial^{2}g_{\gamma\alpha}%
}{\partial x^{\beta}\partial x^{n}}+\frac{\partial^{2}g_{\gamma\beta}%
}{\partial x^{\alpha}\partial x^{n}}-\frac{\partial^{2}g_{\alpha\beta}%
}{\partial x^{\gamma}\partial x^{n}}\right) \\
&  =-g^{\delta\gamma}\left(  \left(  L_{\gamma\alpha}\right)  _{\beta}+\left(
L_{\gamma\beta}\right)  _{\alpha}-\left(  L_{\alpha\beta}\right)  _{\gamma
}\right) \\
&  =0
\end{align*}
where $L_{\alpha\beta}$ is the second fundamental form of the boundary
$\partial M$ and $L_{\alpha\beta}=0$ since $\partial M$ is totally geodesic.
In the same way, we have $\frac{\partial}{\partial x^{n}}\overline{\Gamma
}_{\alpha\beta}^{\delta}\left(  \overline{g}\right)  =0$ on $\partial M$. Then
$\overline{W}_{\alpha\beta n}\left(  x_{0}\right)  =\frac{1-t}{n-2}\left(
\overline{u}_{nnn}+v_{nnn}\right)  \overline{g}_{\alpha\beta}$. Similarly, by
Lemma 5 we have
\[
\overline{W}_{nnn}\left(  x_{0}\right)  =v_{nnn}\left(  x_{0}\right)
+\frac{1-t}{n-2}\left(  \overline{u}_{nnn}+v_{nnn}\right)  \left(
x_{0}\right)  .
\]
Now differentiating (\ref{4.17}) alone the normal direction and taking its
value at $x_{0}$ we have
\begin{align*}
e^{2\overline{u}\left(  x_{0}\right)  }f_{n}\left(  x_{0},\overline{u}\left(
x_{0}\right)  +v\left(  x_{0}\right)  \right)   &  =F^{\alpha\beta}%
\overline{W}_{\alpha\beta n}\left(  x_{0}\right)  +F^{nn}\overline{W}%
_{nnn}\left(  x_{0}\right) \\
&  =F^{nn}v_{nnn}\left(  x_{0}\right)  +\frac{1-t}{n-2}\left(  \overline
{u}_{nnn}+v_{nnn}\right)  \left(  \sum_{i=1}^{n}F^{ii}\right)  \left(
x_{0}\right)  ,
\end{align*}
where we have used the fact that $F^{\alpha n}(x_{0})=0$. Without loss of
generality, one may assume $u_{nnn}=\overline{u}_{nnn}+v_{nnn}\leq 0$. Then by
use of the condition (A4) and $|\nabla f|\leq\Lambda f$, we have%
\begin{align*}
v_{nnn}\left(  x_{0}\right)   &  \geq e^{2\overline{u}\left(  x_{0}\right)  }%
\frac{f_{n}\left(  x_{0},\overline{u}\left(  x_{0}\right)  +v\left(
x_{0}\right)  \right)  }{F^{nn}}\\
&  \geq -\frac{\Lambda}{\varepsilon}\sigma_{1}\left(  \overline{W}\right)
\geq-C\left(  K+1\right)
\end{align*}
where $C$ depends only on the constants $\Lambda,\varepsilon$, and $a,b$.
Although the covariate derivative is taken with respect to the metric
$\overline{g}$, it is the same if we take the covariate derivative with
respect to the metric $g$ on the boundary $\partial M$, i.e. $v_{nnn}%
(g)=v_{nnn}(\overline{g})$. Now we have $u_{nnn}\left(  x_{0}\right)
=\overline{u}_{nnn}\left(  x_{0}\right)  +v_{nnn}\left(  x_{0}\right)
\geq-C\left(  K+1\right)  $.\begin{flushright}$\Box$
\end{flushright}\textit{Step 2. }

By Step 1, we have shown that the maximum point of $\overline{H}$ must be in
the interior of $M$. Then similar with the computation of Theorem 2, we have
at the maximam point $x_{0}$
\[
0=\overline{H}_{i}=e^{px_{n}}(\eta_{i}K+\eta K_{i}+p\delta_{in}K\eta)
\]
That is
\[
K_{i}=-(\frac{\eta_{i}}{\eta}+p\delta_{in})K.
\]
We also have
\begin{align*}
0  &  \geq\overline{H}_{ij}=e^{px_{n}}((\eta_{i}K+\eta K_{i}+p\delta_{in}%
K\eta)p\delta_{jn}+\eta_{ij}K\\
&  +\eta K_{ij}+\eta_{i}K_{j}+\eta_{j}K_{i}+p\delta_{in}K_{j}\eta+p\delta
_{in}K\eta_{j}).
\end{align*}
Note that $|\nabla\eta|\leq\frac{C\eta^{1/2}}{r},|\nabla^{2}\eta|\leq\frac
{C}{r^{2}}$, we have
\[
0\geq\overline{H}_{ij}=e^{px_{n}}(\eta K_{ij}+\overline{\Lambda}_{ij}K),
\]
where
\begin{align*}
\overline{\Lambda}_{ij}  &  =\eta_{ij}-p\eta_{i}\delta_{jn}-p\eta_{j}%
\delta_{in}-p^{2}\eta\delta_{in}\delta_{jn}-\frac{2\eta_{i}\eta_{j}}{\eta}\\
&  \geq-C\left(  p^{2}+1\right)  \delta_{ij}%
\end{align*}
and $C$ depends only on $r$. By the above inequalities, similar with
(\ref{4.4}) we have%
\begin{align}
0  &  \geq\eta P^{ij}\overline{H}_{ij}e^{-px_{n}}\nonumber\\
&  =\eta^{2}P^{ij}K_{ij}+\eta\overline{\Lambda}_{ij}P^{ij}K\nonumber\\
&  \geq\eta^{2}P^{ij}\sum_{k}[u_{ijkk}+2a(u_{ki}u_{kj}+u_{ijk}u_{k}%
)]-C(\sum_{i}F^{ii})(1+\left\vert \nabla^{2}u\right\vert ^{3/2}).
\tag{4.18}\label{4.18}%
\end{align}
We estimate the terms $\sum_{k}P^{ij}u_{ijkk}$ and $\sum_{k}P^{ij}%
(u_{ki}u_{kj}+u_{ijk}u_{k})$ respectively. As the proof of Theorem 2 (b), we
may get (\ref{4.13}) and (\ref{4.14}). Then by the cancellations, we may get
(\ref{4.15}). Therefore we get the estimations of $\left\vert \nabla
^{2}u\right\vert $ and $|\nabla u|^{2}$. \begin{flushright}$\Box$
\end{flushright}

\noindent
\section{Proof of Theorems 3 and 4}

\noindent Similar as the proofs of Theorems 1 and 2, we first prove Theorem 4
for some general functions $a\left(  x\right)  $ and $b\left(  x\right)  $,
and a general $2$-symmetric tensor $S$. We then study the boundary estimates
for special functions $a\left(  x\right)  $, $b\left(  x\right)  $ and a
special $2$-symmetric tensor $S$, i.e. $a\left(  x\right)  ,b\left(  x\right)
$ are both constants and $S=A$ is the Schouten tensor.

\noindent\textit{Proof of Theorem 4. }

\noindent\textit{(1) Case (a).}

\noindent Let $H=\eta(\triangle u+a|\nabla u|^{2})$ and $K=\triangle
u+a|\nabla u|^{2}$, where $0\leq\eta\leq1$ is a cutoff function as before.

Note that $\Gamma\subset\Gamma_{1}^{+}$ and $V=\frac{t-1}{n-2}(\triangle
u)\ g-\nabla^{2}u-a(x)du\otimes du-b(x)|\nabla u|^{2}g+S,$ we can immediately
get
\[
0\leq\ tr(V)=(n\frac{t-1}{n-2}-1)\triangle u-(a+nb)|\nabla u|^{2}%
+trS\leq(n\frac{t-1}{n-2}-1)K-n\delta_{1}|\nabla u|^{2}+C,
\]
Hence, $|\nabla u|^{2}\leq\frac{(n\frac{t-1}{n-2}-1)K+C}{n\delta_{1}}$. Thus
we have
\begin{equation}
|\nabla u|^{2}\leq C(K+1), \tag{5.1}\label{5.1}%
\end{equation}
where $C$ depends only on $||a||_{\infty},||b||_{\infty}$ and $\delta_{1}$. By
(\ref{5.1}), we can obtain $\triangle u<C(K+1)$ and
\begin{equation}
|\nabla^{2}u|\leq C(K+1). \tag{5.2}\label{5.2}%
\end{equation}
Let $x_{0}$ be an interior point where $H$ achieves its maximum. At $x_{0}$,
we have
\[
0=H_{i}=\eta_{i}K+\eta K_{i},
\]
that is
\[
K_{i}=-\frac{\eta_{i}}{\eta}K.
\]
We also have%
\[
0\geq H_{ij}=\eta_{ij}K+\eta K_{ij}+\eta_{i}K_{j}+\eta_{j}K_{i}.
\]
Note that $|\nabla\eta|\leq\frac{C\eta^{1/2}}{r},|\nabla^{2}\eta|\leq\frac
{C}{r^{2}}$, we have
\[
0\geq H_{ij}\triangleq\eta K_{ij}+\Lambda_{ij}K,
\]
where $\Lambda_{ij}$ is bounded. If we take
\[
Q^{ij}=\frac{t-1}{n-2}(\sum F^{ll})\delta_{ij}-F^{ij},
\]
which is also positive definite when $t\ge n-1$, we can obtain
\[
0\geq\eta Q^{ij}H_{ij}=-\eta F^{ij}H_{ij}+\eta\frac{t-1}{n-2}(\sum
F^{ii})H_{kk}.
\]
By the same computation as in the case (a) of Theorem 1, we may get
\begin{equation}
0\geq\sum F^{ii}(2(\frac{t-1}{n-2}a+b)\eta^{2}|\nabla^{2}u|^{2}-C\eta
^{3/2}K^{3/2}-C\eta K-C). \tag{5.3}\label{5.3}%
\end{equation}
As in the case (a) of Theorem 2, we may discuss (\ref{5.3}) in two cases. If
there exists a constant $A>0$, such that $|\nabla u|^{2}({x_{0}})<A|\triangle
u|({x_{0}})$, we may prove\textbf{\ }%
\[
|\triangle u|({x_{0}})\leq C
\]
and
\[
K\leq C.
\]
Otherwise for any constant $A>0$ large enough, $|\nabla u|^{2}({x_{0}})\geq
A|\triangle u|({x_{0}})$. By use of the assumption that $\min\left\{
2ab+b^{2},b^{2}\right\}  \geq\delta_{2}>0$, we may prove
\[
|\nabla u|^{2}({x_{0}})\leq C,
\]
therefore we have
\[
K\leq C.
\]
By (\ref{5.2}), we get the Hessian estimates. \bigskip

\noindent\textit{(2) Case (b).}

We take the same auxiliary function $H=\eta(\triangle u+a|\nabla
u|^{2})\triangleq\eta K$ as in the case (a), where $0\leq\eta\leq1$ is a
cutoff function such that $\eta=1$ in $B_{\frac{r}{2}}$ and $\eta=0$ outside
$B_{r}$, and also $|\nabla\eta|\leq\frac{C\eta^{1/2}}{r},|\nabla^{2}\eta
|\leq\frac{C}{r^{2}}$.

Since $a(x)+nb(x)\geq\delta_{3}$, by the condition $\Gamma\subset\Gamma
_{1}^{+}$ again, we have
\begin{align*}
0  &  \leq\ tr(V)=(n\frac{t-1}{n-2}-1)\triangle u-a|\nabla u|^{2}-nb|\nabla
u|^{2}+trS\\
&  \leq(n\frac{t-1}{n-2}-1)\triangle u-\delta_{3}|\nabla u|^{2}+C,
\end{align*}
and then
\begin{equation}
|\nabla u|^{2}\leq C(\triangle u+1). \tag{5.4}\label{5.4}%
\end{equation}
Without loss of generality, we may assume
\[
K=\Delta u+a\left\vert \nabla u\right\vert ^{2}>>1.
\]
Since $a\left(  x\right)  \geq0$, by (\ref{5.4}), we have \
\begin{equation}
\Delta u\leq C(K+1) \tag{5.5}\label{5.5}%
\end{equation}
and
\begin{equation}
|\nabla u|^{2}\leq C(K+1). \tag{5.6}\label{5.6}%
\end{equation}
Suppose that the maximum point of $H$ achieves at $x_{0}$, an interior point,
we may get an inequality just replacing $K$ in (\ref{5.3}) by $\left\vert
\nabla^{2}u\right\vert $
\begin{equation}
0\geq\sum F^{ii}(2(\frac{t-1}{n-2}a+b)\eta^{2}|\nabla^{2}u|^{2}-C(\eta
|\nabla^{2}u|)^{3/2}-C(\eta|\nabla^{2}u|)-C). \tag{5.7}\label{5.7}%
\end{equation}
The coefficient of the highest order term $\frac{t-1}{n-2}a(x)+b(x)\ge \frac
{\delta_{3}}{n}>0$ since $a(x)\geq0$ and $a(x)+nb(x)\ge\delta_{3}>0$. Therefore
we can get the bounds of $K,\left\vert \nabla^{2}u\right\vert $ and $|\nabla
u|^{2}$. \begin{flushright} $\Box$
\end{flushright}

\noindent\textit{Proof of Theorem 3.}

Note that $a,b$ are two constants, $-S=A$ is the Schouten. Similar as the
proof of Theorem 1 Case (b), by (\ref{5.5}) and (\ref{5.6}), we only need to
estimate $K=\triangle u+a|\nabla u|^{2}$. Consider $\overline{H}=\eta
Ke^{px_{n}}$, where $0\leq\eta\leq1$ is a cutoff function as before. We may
show the maximum point of $\overline{H}$ must be in the interior of $M$. Then
the argument in Theorem 4 case (b) to get the estimations.

We prove this by contradiction. Assume the maximum point of $\overline{H}$,
$x_{0}$, is on the boundary, then by (\ref{2.1}), (\ref{2.2}) and (\ref{2.4}),
we have $u_{\alpha\alpha n}+2au_{\alpha n}u_{\alpha}+2au_{nn}u_{n}%
+a_{n}(u_{\gamma}u_{\gamma}+u_{n}u_{n})=0.$ Then
\[
\overline{H}_{n}|_{x_{0}}=\eta e^{px_{n}}\left(  u_{nnn}+pK\right)  .
\]
Furthermore, we can get the following Lemma 7 as well:

\noindent\textbf{Lemma 7. }\ \textit{We can find some positive constant $C$,
such that $u_{nnn}(x_{0})\geq-C(K+1)$. }

From Lemma 7, we can show that $\overline{H}_{n}|_{x_{0}}>0$ as long as $p$ is
large enough, which contradicts with the assumption that $x_{0}$ is a maximum
point. Hence, $\overline{H}$ achieves its maximum at an interior point.

\noindent\textit{Proof of Lemma 7.}\textbf{ }

We may assume $u_{nnn}\leq0$. Similar as Lemma 6, by Lemma 5, we may choose a
conformal metric $\bar{g}=e^{-2\bar{u}}g$ and $\overline{u}_{n}|_{\partial
M}=0$ at first. In this metric, $\partial M$ is still totally geodesic and
$\overline{A}_{\alpha\beta,n}(x_{0})=0$. We wish to find a metric
$\widetilde{g}=e^{-2v}\overline{g}$ such that $u=\overline{u}+v$ is a solution
to (\ref{1.7}). Then equation (\ref{1.7}) becomes
\begin{equation}
\left\{
\begin{array}
[c]{ll}%
F(\bar{g}^{-1}\overline{V})=e^{2\overline{u}}f(x,\overline{u}+v) & \text{
in\ }\overline{B}_{r}^{+},\\
\frac{\partial v}{\partial x^{n}}=0 & \text{ on }\Sigma_{r},
\end{array}
\right.  \tag{5.8}\label{5.8}%
\end{equation}
where
\begin{align*}
\overline{V}_{ij}  &  =\frac{t-1}{n-2}\left(  \overline{\Delta}\left(
\overline{u}+v\right)  +\overline{g}^{lk}\left(  \overline{\Gamma}_{lk}%
^{p}\left(  \overline{g}\right)  -\Gamma_{lk}^{p}\left(  g\right)  \right)
\left(  \overline{u}_{p}+v_{p}\right)  \right)  \overline{g}_{ij}\\
&  -\left(  \overline{A}_{ij}+\left(  a-1\right)  \overline{u}_{i}\overline
{u}_{j}+\left(  b+\frac{1}{2}\right)  \left\vert \overline{\nabla}\overline
{u}\right\vert _{\overline{g}}^{2}\overline{g}_{ij}\right) \\
&  -\left(  \overline{\nabla}_{ij}^{2}v+\left(  \overline{\Gamma}_{ij}%
^{k}\left(  \overline{g}\right)  -\Gamma_{ij}^{k}\left(  g\right)  \right)
v_{k}+a\overline{u}_{i}v_{j}+a\overline{u}_{j}v_{i}+av_{i}v_{j}\right) \\
&  -\left(  b\left(  \overline{u}_{k}v_{l}+\overline{u}_{l}v_{k}\right)
\overline{g}^{kl}\overline{g}_{ij}+b\left\vert \overline{\nabla}v\right\vert
_{\overline{g}}^{2}\overline{g}_{ij}\right)  .
\end{align*}
Notice that the boundary $\partial M$ preserves totally geodesic, we have
$u_{n}=\overline{u}_{n}=0$, and $u_{n\alpha}=\overline{u}_{n\alpha}=0$. By
Lemma 3, we have $u_{\alpha\beta n}=\overline{u}_{a\beta n}=0$, therefore
$v_{n}=v_{n\alpha}=v_{\alpha\beta n}=0$ on $\partial M$. As lemma 6, we have
$\overline{V}_{\alpha n}(x_{0})=0$. Employing an argument of Lemma 13 in
\cite{Cn3}, we know $F^{\alpha n}(x_{0})=0$. By Lemma 5, similar as the
computation in the proof of Lemma 6, we have
\[
\overline{V}_{\alpha\beta n}\left(  x_{0}\right)  =\frac{t-1}{n-2}\left(
\overline{u}_{nnn}+v_{nnn}\right)  \overline{g}_{\alpha\beta}%
\]
and
\[
\overline{V}_{nnn}\left(  x_{0}\right)  =-v_{nnn}\left(  x_{0}\right)
+\frac{t-1}{n-2}\left(  \overline{u}_{nnn}+v_{nnn}\right)  \left(
x_{0}\right)  .
\]
Differentiating (\ref{5.8}) alone the normal direction and taking its value at
$x_{0}$ we have
\begin{align*}
e^{2\overline{u}\left(  x_{0}\right)  }f_{n}\left(  x_{0},\overline{u}\left(
x_{0}\right)  +v\left(  x_{0}\right)  \right)   &  =F^{\alpha\beta}%
\overline{V}_{\alpha\beta n}\left(  x_{0}\right)  +F^{nn}\overline{V}%
_{nnn}\left(  x_{0}\right) \\
&  =\frac{t-\left(  n-1\right)  }{n-2}F^{nn}v_{nnn}\left(  x_{0}\right)
+\frac{t-1}{n-2}F^{nn}\overline{u}_{nnn}\left(  x_{0}\right) \\
&  +\frac{t-1}{n-2}\left(  \overline{u}_{nnn}+v_{nnn}\right)  \left(
\sum_{\alpha=1}^{n-1}F^{\alpha\alpha}\right)  \left(  x_{0}\right)  .
\end{align*}
Since we have assumed that $u_{nnn}\left(  x_{0}\right)  \leq0$, this means
that $\left(  \overline{u}_{nnn}+v_{nnn}\right)  \left(  x_{0}\right)  \leq0.$
We therefore have
\begin{align*}
v_{nnn}\left(  x_{0}\right)   &  \geq\frac{n-2}{t-\left(  n-1\right)  }\left[
e^{2\overline{u}\left(  x_{0}\right)  }\frac{f_{n}\left(  x_{0},\overline
{u}\left(  x_{0}\right)  +v\left(  x_{0}\right)  \right)  }{F^{nn}}-\frac
{t-1}{n-2}\overline{u}_{nnn}\left(  x_{0}\right)  \right] \\
&  \geq-C\left(  K+1\right)  ,
\end{align*}
since $t>n-1$, where we have used the condition (A4) and $|\nabla f|\leq\Lambda f$, the
constant $C$ depends only on the constants $\Lambda,\varepsilon$, and $a,b,t$.
\begin{flushright} $\Box$
\end{flushright}

\bigskip Now $\overline{H}=\eta Ke^{px_{n}}$ attains its maximum at an
interior point $x_{0}$, we have at $x_{0}$
\[
0=\overline{H}_{i}=e^{px_{n}}(\eta_{i}K+\eta K_{i}+p\delta_{in}K\eta),
\]
that is
\[
K_{i}=-(\frac{\eta_{i}}{\eta}+p\delta_{in})K.
\]
We also have%
\begin{align*}
0  &  \geq\overline{H}_{ij}=e^{px_{n}}((\eta_{i}K+\eta K_{i}+p\delta_{in}%
K\eta)p\delta_{jn}+\eta_{ij}K\\
&  +\eta K_{ij}+\eta_{i}K_{j}+\eta_{j}K_{i}+p\delta_{in}K_{j}\eta+p\delta
_{in}K\eta_{j}).
\end{align*}
Then
\[
0\geq\overline{H}_{ij}\triangleq e^{px_{n}}(\eta K_{ij}+\overline{\Lambda
}_{ij}K),
\]
where $\overline{\Lambda}_{ij}$ is bounded. Taking
\[
Q^{ij}=\frac{t-1}{n-2}(\sum F^{ll})\delta_{ij}-F^{ij},
\]
as the proof of Theorem 4, we can obtain
\[
0\geq\eta Q^{ij}\overline{H}_{ij}e^{-px_{n}}=-\eta F^{ij}\overline{H}%
_{ij}e^{-px_{n}}+\eta\frac{t-1}{n-2}(\sum F^{ii})\overline{H}_{kk}e^{-px_{n}%
}.
\]
By the same argument as Case (b) of Theorem 4, we have (\ref{5.7}). Therefore
we get the estimations of $\left\vert \nabla^{2}u\right\vert $ and $|\nabla
u|^{2}$. \begin{flushright} $\Box$
\end{flushright}

\noindent\textit{Addresses}:

\noindent{\small Yan He: Centre for Mathematical Sciences, Zhejiang
University, Hangzhou 310027, China. }

{\small \vskip2pt }

{\small \noindent Weimin Sheng: Department of Mathematics, Zhejiang
University, Hangzhou 310027, China; and Centre for Mathematics and its
Applications, the Australian National University, Canberra, ACT 0200,
Australia.} {\small \vskip5pt }

{\small \noindent\textit{E-mail}: helenaig@hotmail.com, weimins@zju.edu.cn}


\begin{thebibliography}{9999}                                                                                             %


\bibitem[Au]{Au}T. Aubin, Equations diff\'{e}rentielles non lin\'{e}aires et
probl\`{m}e de Yamabe concernant la courbure scalaire, J. Math. Pures Appl.
(9) 55 (1976), 269--296, MR0431287, Zbl0336.53033.

\bibitem[CNS]{CNS}L.A. Caffarelli, L. Nirenberg, and J. Spruck, Dirichlet
problem for nonlinear second order elliptic equations III, Functions of the
eigenvalues of the Hessian, Acta Math. \textbf{1985}, \textit{155}, 261--301,
MR0806416, Zbl 0654.35031.

\bibitem[CGY1]{CGY1}A. Chang, M. Gursky, P. Yang, An equation of Monge-Am
p\'{e}re type in conformal geometry, and four-manifolds of positive Ricci
curvature, Ann. of Math. (2) 155(2002), 709--787, MR1923964, Zbl 1031.53062.

\bibitem[CGY2]{CGY2}A. Chang, M. Gursky, P. Yang, An a priori estimate for a
fully nonlinear equation on four-manifolds, J. Anal. Math. 87 (2002),
151--186, MR1945280, Zbl 1067.58028.

\bibitem[Cn1]{Cn1}S. Chen, Local estimates for some fully nonlinear elliptic
equations, Int. Math. Res. Not. 2005, no. 55, 3403--3425, MR2204639, Zbl pre05017507.

\bibitem[Cn2]{Cn2}S. Chen, Boundary value problems for some fully nonlinear
elliptic equations. \textit{Calc. Var. Partial Differential Equations}, 2007,
\textbf{30(1)}:1--15.

\bibitem[Cn3]{Cn3}S. Chen, Conformal Deformation on Manifolds with Boundary.
Geom. Funct. Anal., 19 (2009), no. 4, 1029--1064.

\bibitem[E1]{E1}Jos\'e F. Escobar, The Yamabe problem on manifolds with
boundary, J. Differential Geom. 35 (1992), no. 1, 21--84.

\bibitem[E2]{E2}Jos\'e F. Escobar, Conformal deformation of a Riemannian
metric to a scalar flat metric with constant mean curvature on the boundary,
Annals of Math., 136 (1992), 1-50.

\bibitem[GeW]{GeW}Y. Ge and G. Wang, On a fully nonlinear Yamabe problem, Ann.
Sci. Ecole Norm. Sup. (4) 39 (2006), 569--598, MR2290138, Zbl pre05125020.

\bibitem[G]{G}B. Guan, Conformal metrics with prescribed curvature curvature
functions on manifolds with boundary, Amer. J. Math., 129 (2007), no. 4, 915--942.

\bibitem[GLW]{GLW}P. Guan, C.-S. Lin and G. Wang, Application of the method of
moving planes to conformally invariant equations. \textit{Math. Z.}, 2004,
\textbf{247(1)}:1--19.

\bibitem[GW1]{GW1}P. Guan and G. Wang, Local estimates for a class of fully
nonlinear equations arising from conformal geometry, Int. Math. Res. Not.
(2003), 1413--1432, MR1976045, Zbl 1042.53021.

\bibitem[GW2]{GW2}P. Guan and G. Wang, A fully nonlinear conformal flow on
locally conformally flat manifolds, J. Reine Angew. Math. 557 (2003),
219--238, MR1978409, Zbl 1033.53058.

\bibitem[GW3]{GW3}P. Guan and G. Wang, Geometric inequalities on locally
conformally flat manifolds, Duke Math. J. 124 (2004), 177--212, MR2072215, Zbl 1059.53034.

\bibitem[GW4]{GW4}P. Guan and X.-J. Wang, On a Monge-Amp\`{e}re equation
arising in geometric optics, J. Diff. Geom., 48(1998), 205--223, MR2072215,
Zbl 0979.35052.

\bibitem[GV1]{GV1}M. Gursky and J. Viaclovsky, A fully nonlinear equation on
four-manifolds with positive scalar curvature, J. Differential Geom. 63
(2003), 131--154, MR2015262, Zbl 1070.53018.

\bibitem[GV2]{GV2}M. Gursky and J.Viaclovsky, Prescribing symmetric functions
of the eigenvalues of the Ricci tensor, Ann. of Math., 2007, \textbf{166}: 475--531.

\bibitem[GV3]{GV3}M. Gursky and J.Viaclovsky, Fully nonlinear equations on
Riemannian manifolds with negative curvature, Indiana Univ. Math. J., 52(3)
(2003), 399-420.

\bibitem[HS1]{HS1}Y. He and W.M. Sheng, On existence of the prescribing $k$
-curvature problem on manifolds with boundary,
Communications in Analysis and Geometry, 19 (2011), no. 1, 53-77.

\bibitem[HS2]{HS2}Y. He and W.M. Sheng, Prescribing the symmetric function of
the eigenvalues of the Schouten tensor, Proc. Amer. Math. Soc., 139 (2011), 1127-1136.

\bibitem[J]{J} Q. Jin, Local Hessian estimates for some conformally invariant fully nonlinear equations with boundary conditions,
Differential and Integral Equations, 20 (2007), no. 2, 121-132.


\bibitem[JLL]{JLL}Q. Jin, A. Li and Y.Y. Li, Estimates and existence results
for a fully nonlinear Yamabe problem on manifolds with boundary, Calc. Var.,
28 (2007), 509-543.

\bibitem[LP]{LP}J. Lee and T. Parker, The Yamabe problem, Bull. Amer. Math.
Soc., 17 (1987), 37-91.

\bibitem[LL1]{LL1}A. Li and Y.Y. Li, On some conformally invariant fully
nonlinear equations, Comm. Pure Appl. Math. 56 (2003), 1416--1464, MR1988895,
Zbl pre02002141.

\bibitem[LL2]{LL2}A. Li and Y.Y. Li, On some conformally invariant fully
nonlinear equations, part II: Liouville, Harnack, and Yamabe, Acta Math. 195
(2005), 117--154, MR2233687, Zbl pre05039005.

\bibitem[LL3]{LL3}A. Li and Y.Y. Li, A fully nonlinear version of the Yamabe
problem on manifolds with boundary, J. Eur. Math. Soc., 8 (2006), 295-316.

\bibitem[LS]{LS}J.Y. Li and W.M. Sheng, Deforming metrics with negative
curvature by a fully nonlinear flow. Calc. Var. PDE., 2005, \textbf{23}: 33--50.

\bibitem[Li]{Li}Y.Y. Li, Local gradient estimates of solutions to some
conformally invariant fully nonlinear equations, C. R. Math. Acad. Sci. Paris
343 (2006), 249--252, MR2245387, Zbl 1108.35061.

\bibitem[LT]{LT}M. Lin and N.S. Trudinger, On some inequalities for elementary
symmetric functions, Bull. Aust. Math. Soc., 50(1994), 317--326, MR1296759,
Zbl 0855.26006.

\bibitem[LTU]{LTU}P.-L. Lions, N.S. Trudinger and J. Urbas, The Neumann
problem for equations of Monge-Amp\'ere type, Comm. Pure Appl. Math., 39
(1986), 539-563.

\bibitem[M]{M}F.C. Marques, A priori estimates for the Yamabe problem in
non-locally conformal flat cases, J. Differential Geom. 71 (2005), no. 2, 315--346.

\bibitem[S1]{S1}R. Schoen, Conformal deformation of a Riemannian metric to
constant scalar curvature, J. Diff. Geom. 20(1984), 479--495, MR0788292, Zbl 0576.53028.

\bibitem[S2]{S2}R. Schoen, Variational theory for the total scalar curvature
functional for Riemannian metrics and related topics. Topics in Calculus of
Variations, Lectures Notes in Math. 1365, Springer,1989, pp. 120--154,
MR0994021, Zbl 0702.49038.

\bibitem[S]{S}W.M. Sheng, Admissible metrics in $\sigma_k$-Yamabe
equation, Proc. Amer. Math. Soc., 2007, \textbf{136(5)}:1795--1802.

\bibitem[STW]{STW}W.M. Sheng, N.S. Trudinger and X-J. Wang, The Yamabe problem
for higher order curvatures, J. Diff. Geom. 77(2007), 515-553.

\bibitem[SZ]{SZ}W.M. Sheng and Y. Zhang. A class of fully nonlinear equations
arising from conformal geometry. Math. Z., 2007, \textbf{\ 255}: 17--34.

\bibitem[Tr1]{Tr1}N.S. Trudinger, Remarks concerning the conformal deformation
of Riemannian structures on compact manifolds, Ann. Scuola Norm. Sup. Pisa (3)
22(1968), 265--274, MR0240748, Zbl 0159.23801.

\bibitem[Tr2]{Tr2}N.S. Trudinger, The Dirichlet problem for the prescribed
curvature equations, Arch. Rational Mech. Anal., 111(2) (1990), 153-179.

\bibitem[TW1]{TW1}N.S. Trudinger and X-J. Wang, On Harnack inequalities and
singularities of admissible metrics in the Yamabe problem, Calc. Var. Partial
Differential Equations, 35 (2009), no. 3, 317-338.

\bibitem[TW2]{TW2}N.S. Trudinger and X-J. Wang, The intermediate case of the
Yamabe problem for higher order curvatures, International Mathematics Research
Notices, 2009.

\bibitem[U]{U} J.Urbas, An expansion of convex hypersurfaces,
J. Differential Geometry, 33(1991), 91-125.

\bibitem[V1]{V1}J. Viaclovsky, Conformal geometry, contact geometry, and the
calculus of variations, Duke Math. J. 101 (2000), 283--316, MR1738176, Zbl 0990.53035.

\bibitem[V2]{V2}J. Viaclovsky, Estimates and existence results for some fully
nonlinear elliptic equations on Riemannian manifolds, Comm. Anal. Geom. 10
(2002), 815--846, MR1925503, Zbl 1023.58021.

\bibitem[W1]{W1}X.-J. Wang, A class of fully nonlinear elliptic equations and
related functionals, Indiana Univ. Math. J., 43 (1994), 25--54, MR1275451, Zbl 0805.35036.

\bibitem[W2]{W2}X.-J. Wang, A priori estimates and existence for a class of
fully nonlinear elliptic equations in conformal geometry, Chinese Ann. Math.
Ser. B, 27 (2006), 169--178, MR2243678, Zbl 1104.53035.

\bibitem[Ya]{Ya}H. Yamabe, On a deformation of Riemannian structures on
compact manifolds, Osaka Math. J. 12(1960), 21--37, MR0125546, Zbl 0096.37201.

\bibitem[Ye]{Ye}R. Ye, Global existence and convergence of Yamabe flow, J.
Differential Geom. 39 (1994), 35--50, MR1258912, Zbl 0846.53027.
\end{thebibliography}
\end{document}